\documentclass{article}
\usepackage{graphicx} % Required for inserting images
\usepackage[margin=1in]{geometry}

\usepackage{hyperref}
\usepackage{amsmath,amssymb,bm}
\usepackage{xcolor}
\usepackage[multiple]{footmisc}

\usepackage{amsthm}

\newcommand{\trans}{\mathsf{T}}
\newcommand{\Rh}{R_h}
\newcommand{\sympmap}{F}
\newcommand{\Tnorm}[1]{\left|#1 \right|_{\mathbb{T}}}

\newcommand{\abs}[1]{\left|#1\right|}
\newcommand{\norm}[1]{\left\lVert#1\right\rVert}

\newcommand{\floor}[1]{\left\lfloor#1\right\rfloor}
\newcommand{\ip}[2]{\left \langle #1 , \, #2 \right \rangle}
\newcommand{\JMAP}{J_0}

\newcommand{\pd}[2]{\frac{\partial #1}{\partial #2}}
\newcommand{\Cbb}{\mathbb{C}}
\newcommand{\Rbb}{\mathbb{R}}
\newcommand{\Tbb}{\mathbb{T}}
\newcommand{\Zbb}{\mathbb{Z}}
\newcommand{\Ocal}{\mathcal{O}}
\newcommand{\I}{i}

\newcommand{\dif}{\, \mathrm{d}}

\DeclareMathOperator{\tr}{tr}
\DeclareMathOperator{\Diag}{Diag}

\DeclareMathOperator*{\argmin}{arg\,min}
\DeclareMathOperator{\prob}{\mathbb{P}}
\DeclareMathOperator{\expec}{\mathbb{E}}
\DeclareMathOperator{\real}{Re}
\DeclareMathOperator{\imag}{Im}

\title{Robust computation of higher-dimensional invariant tori from individual trajectories}
\author{Maximilian Ruth$^{1,2}$, Jackson Kulik$^3$, Joshua Burby$^4$}
\date{
$^1$Contact: \texttt{maximilian.ruth@austin.utexas.edu}\\
$^2$Oden Institute for Computational Engineering and Sciences, University of Texas at Austin, Austin, TX 78712, USA\\
$^3$Department of Mechanical and Aerospace Engineering, Utah State University, Logan, UT 84322, USA\\
$^4$Department of Physics and Institute for Fusion Studies, The University of Texas at Austin, Austin, TX 78712, USA\\[2ex]
\today}

\begin{document}

\maketitle

\section*{Abstract}
We present a method for computing invariant tori of dimension greater than one.
The method uses a single short trajectory of a dynamical system without any continuation or initial guesses. 
No preferred coordinate system is required, meaning the method is practical for physical systems where the user does not have much \textit{a priori} knowledge.
Three main tools are used to obtain the rotation vector of the invariant torus: the reduced rank extrapolation method, Bayesian maximum a posteriori estimation, and a Korkine-Zolatarev lattice basis reduction.
The parameterization of the torus is found via a least-squares approach.
The robustness of the algorithm is demonstrated by accurately computing many two-dimensional invariant tori of a standard map example.
Examples of islands and three-dimensional invariant tori are shown as well.

\section{Introduction}
Invariant tori are central to the study of low-dimensional Hamiltonian dynamical systems.
% Trajectories initialized on invariant tori are guaranteed predictable bounded motion, which is often appealing in application.
Familiar examples in celestial mechanics are elliptic orbits in the two-body problem and quasiperiodic trojan orbits in the restricted three-body problem \cite{murray_solar_2000}.
In plasma confinement, invariant tori appear for both magnetic fieldline dynamics and charged particle motion \cite{hazeltine2003,imbert-gerard_introduction_2024}.
For both areas, the predictable bounded motion of integrable trajectories is critical for downstream physics and engineering problems.
In this paper, we will parameterize these tori.

There is a variety of existing methods for computing invariant tori.
The most common one is the parameterization method \cite{haro2016,haro_parameterization_a,haro_parameterization_b,haro_parameterization_c} which directly solves the conjugacy that defines an invariant torus.
One of the primary benefits of the parameterization method is its connection to Kolmogorov-Arnold-Moser (KAM) theory \cite{LLave_2004,kolmogorov-conservation-1954,arnold2009proof,moser1962invariant}, where a convergent numerical algorithm amounts to a proof of the existence of a torus \cite{Figueras2017}.
The parameterization method is used within astrodynamics for the computation of quasiperiodic orbits \cite{Baresi2018,kolemen2012, jorba2009computation, mccarthy2021leveraging}
A similar method known as quadratic flux minimizing (QFMin) surfaces \cite{dewar1992,dewar1994} is used for flux surfaces in stellarators \cite{hudson_almost-invariant_1996,hudson_reduction_2001} (see \cite{giuliani_direct_2022} for another approach).
QFMin recovers an invariant torus when it exists and otherwise it returns optimal surfaces that cut through island chains and chaos. 

The appealing aspect of both the parameterization method and QFMin is the accuracy with which an invariant torus can be computed. 
However, these methods also come with a practical downside: a good initial guess is required.
This is typically solved by continuation either in a physical parameter or within phase space.
For instance, given the location of an elliptic fixed point of a 2D symplectic map (an ``o-point''), one can use continuation to compute nested invariant tori about that point.
Alternatively, parameter continuation is typically used to track invariant tori as the magnitude of some forcing term increases. % \cite{greene1979,Figueras2017??}.
In higher dimensional systems, continuation becomes significantly more complicated, as initial guesses are not as readily available.
This is a particular problem when one wants to compute many tori, as continuation often requires a significant amount of effort from the researcher.

Another class of methods---which we refer to as frequency methods---are based on single trajectories of the dynamical system.
The most popular work within this category is the frequency map method by Laskar \cite{laskar1999,laskar_chaotic_1990}, which can be used to compute rotation numbers of invariant tori from individual trajectories.
With a similar formulation, it was shown that the weighted Birkhoff average \cite{Das2017,sander2020,meiss_resonance_2025} can be used to accelerate the computation of $C^\infty$ ergodic averages of length $T$ trajectories on $C^\infty$ tori with a rate of $\Ocal(T^{-M})$ for all $M$. 
This result has been extended in the analytic case to exponential rates of $\Ocal(e^{-r T^\xi})$ in \cite{Tong2024} (see \cite{tong_quantitative_2024} for sharper results). 
Using an ergodic inverse Fourier transform, the weighted Birkhoff average can then be used to compute a parameterization of an invariant torus in two dimensions \cite{blessing_weighted_2024} (see also Section \ref{sec:parameterization}).
This approximation of the invariant torus can be used to initialize the parameterization method for a higher fidelity result.
Here, we will rely upon the Birkhoff reduced rank extrapolation (Birkhoff RRE) method \cite{ruth_finding_2024}, which simultaneously computes many frequencies to high precision rather than the single frequency in frequency map analysis.
Birkhoff RRE has been used to accurately compute magnetic flux surfaces in \cite{ruth_high-order_2025} and charged particle orbits in \cite{burby_nonperturbative_2025}.

Frequency methods have the benefit of using a single trajectory, rather than using continuation as is typical for the parameterization method and QFMin.
Unfortunately, frequency methods typically require the assumption that the researcher has \textit{a priori} knowledge that the invariant torus is star shaped to obtain the rotation number.
For simple one-dimensional invariant circles, this is not an unreasonable assumption; a na\"ive average of the trajectory position will often give a good enough guess for a coordinate system.
In the case that a 1D trajectory is not star-shaped, the method by Luque and Villanueva \cite{luque_numerical_2009} can be used to alleviate the problem.
However, to our knowledge, there is no general frequency method for higher dimensional tori.
In particular, for the rotation vector to be easily expressed via a limit, it is necessary to know a coordinate system where the torus is a graph over angular variables.

In this paper, we present a new frequency method to compute invariant tori that satisfies
\begin{enumerate}
    \item The method is dimension- and coordinate-agnostic.
    \item No continuation or initial guesses are required for the torus parameterization or the rotation vector.
    \item The input to the algorithm is a single trajectory on an invariant torus.
    \item Few iterations of the map are required (as compared to weighted Birkhoff averaging).
\end{enumerate}
We emphasize that none of the reviewed methods satisfy all of these requirements.
The first and second requirements are the most difficult to perform simultaneously, while the third and fourth are convenient computationally.
The fourth point is particularly appealing for many physical systems, as computing long trajectories may be the most costly part of computing tori in practice.
Because we satisfy the above requirements, this method allows for more practical computations of higher-dimensional invariant tori within several application areas (code is publically available in the \texttt{SymplecticMapTools.jl} Julia package \cite{Ruth:2025:SymplecticMapTools}).

Our method comprises four steps.
First, we accurately compute the trajectory's most prominent frequencies using Birkhoff RRE. 
This step has the side benefit of returning a residual that can be used to classify the trajectory as integrable or chaotic.
Second, we use those frequencies to infer a valid rotation vector via Bayesian maximum \textit{a posteriori} (MAP) estimation.
Third, we transform the rotation vector so that the invariant torus has a preferred coordinate system. 
(This aligns with the intuition of a ``short way'' and ``long way'' round a two-dimensional torus.)
This is made precise via a Korkine-Zolatarev (KZ) lattice basis reduction using an averaged metric on the torus.
Fourth, we compute the invariant torus Fourier coefficients from the trajectory and the rotation vector. 

Because the final parameterization step is the simplest, we begin our exposition with it in Section \ref{sec:parameterization} assuming the rotation vector is already known.
This will allow us to explore how the parameterization of higher dimensional tori changes under changes in the homology. 
Then, in Section \ref{sec:rotation-vector} we explain the first three steps to compute the rotation vector. 
In both of these sections, we use examples that demonstrate how high dimensionality introduces difficulties not seen for one-dimensional invariant tori.
We test the method in Section \ref{sec:examples} with three applications: a coupled standard map, charged particle dynamics in a stellarator, and non-planar satellite trajectories in cislunar space.
Through these tests, we explore the robustness of the method in situations relevant to the modern theory of dynamical systems, plasma physics, and astrodynamics. 
By one \textit{a posteriori} measure of success, we show that we successfully parameterize 95\% of the initialized tori for the weakly coupled standard map, and we discuss the two main modes of failure that are accentuated in higher dimensions: nearly resonant rotation numbers and filamentary invariant tori.
Finally, we conclude in Section~\ref{sec:conclusions}.

\section{Computing invariant torus Fourier coefficients}
\label{sec:parameterization}
Let $\sympmap : \mathcal{M} \to \mathcal{M}$ be an smooth map on a smooth manifold $\mathcal{M}$ (here, smooth either means analytic or $C^M$ for some $M$, which we will discuss when relevant).
We define a $d$-dimensional invariant torus as an smooth mapping $S : \Tbb^d \to \mathcal{M}$ where $\Tbb = \Rbb \backslash \Zbb$ such that the dynamics are conjugate to rotation by a vector $\bm \omega \in \Tbb^d$ as
\begin{equation*}
    \forall\bm{\theta}\in \mathbb{T}^d,\quad \sympmap(S(\bm \theta)) = S(\bm \theta + \bm \omega).
\end{equation*}
We assume that $\bm \omega$ obeys a Diophantine condition, where for all $\bm k \in \Zbb^d$
\begin{equation*}
    \Tnorm{\bm \omega \cdot \bm k} \geq \frac{c}{\norm{\bm k}^{\nu}}
\end{equation*}
for some constants $c>0$ and $\nu\geq d$, where the vector norm $\norm{\star}$ is $L^2$ unless otherwise specified, and $\Tnorm{\star} : \Tbb \to \Rbb$ is the distance of an angle $\theta \in \Tbb$ from zero 
\begin{equation*}
    \Tnorm{\theta} = \min_{m\in\Zbb} \abs{\theta - m}.
\end{equation*}
Note that we regard $\bm{\omega}\cdot\bm{k}$ as an element of $\mathbb{T}$. The Diophantine condition guarantees the rotation vector $\bm{\omega}$ is far from being resonant ($\bm \omega$ is resonant if $\bm \omega \cdot \bm k = 0$ for some $\bm k \in \Zbb^d\backslash\{0\}$).
Non-resonant $\bm{\omega}$ implies ergodic dynamics on the invariant torus. 
Invariant tori commonly arise in Poincar\'e maps for Hamiltonian dynamical systems. 
For this reason, we focus on symplectic maps within this document, as continuous-time tori can be found by simply ``extruding'' the discrete-time tori via the flow.

An example of a $d=2$ invariant torus is found in the elliptic restricted three-body problem (ER3BP) \cite{szebehely_elliptic_1964}.
The ER3BP describes the motion of a body of negligible mass forced by the gravity of two massive primary bodies in an elliptic two-body orbit about one another. 
In our example, we choose the primary body parameters to match Jupiter and the Sun.
The frame of reference is chosen to co-rotate and compress with the elliptic orbit, giving the pulsating synodic Cartesian coordinates $\bm \xi = (\xi,\eta)$ in the plane of the orbit.
In these coordinates, the dynamics of the third body are given by the 2.5-dimensional Hamiltonian dynamical system
\begin{equation}
\label{eq:rer3bp}
    \dot{\bm \xi} = -\pd{H}{\bm p}, \qquad \dot{\bm p} = \pd{H}{\bm \xi}, \qquad H(\bm \xi, \bm p, f) = \frac{1}{2}\left(p_\xi^2 + p_\eta^2 + \xi^2 + \eta^2 \right) + \eta p_\xi - \xi p_\eta - \phi(\bm \xi,f),
\end{equation}
where $\bm p = (p_\xi,p_\eta)$ are the momenta corresponding to $\bm \xi$ and the dot represents derivatives with respect to the time-like true anomaly $f$ --- defined as the angle between the periapsis of either of the primary bodies and its current position relative to the common center-of-mass of the system.
The potential $\phi$ is defined by 
\begin{equation*}
    \phi = \frac{1}{1+ \epsilon \cos f}\left[\frac{1}{2}\left((1-\mu)\rho_1^2 + \mu \rho_2^2\right) + \frac{1-\mu}{\rho_1} + \frac{\mu}{\rho_2} \right], \quad \rho_1^2 = (\xi-\mu)^2 + \eta^2, \quad \rho_2^2 = (\xi - (\mu-1))^2 + \eta^2,
\end{equation*}
where $\epsilon = 0.0489$ is the eccentricity of the elliptical orbit of the primary bodies and $\mu = m_2/(m_1+m_2) = 9.54\times 10^{-4}$ is the mass ratio of the primaries.
A map is obtained by evolving the ER3BP over an orbit of the two primaries, i.e.~$\sympmap : (\bm \xi(0), \bm p(0)) \mapsto (\bm \xi(2\pi), \bm p(2\pi))$.
To compute these trajectories, we use the SRK3 symplectic integrator \cite{zhao_unified_2014} provided by the Julia package GeometricIntegrators.jl \cite{Kraus:2020:GeometricIntegrators} with a time step of $2\pi/1000$.

\begin{figure}
    \centering
    \includegraphics[width=0.9\linewidth]{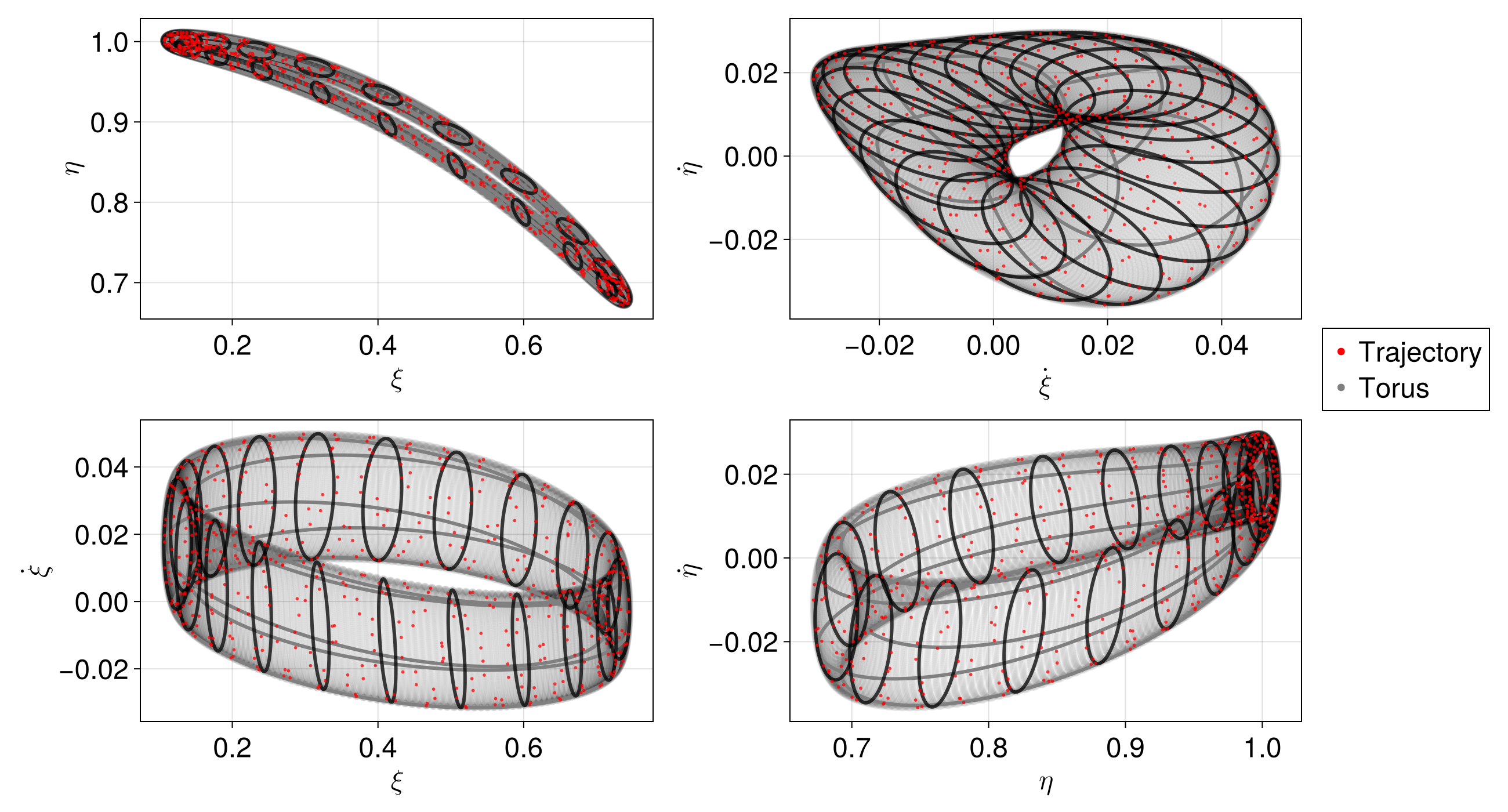}
    \caption{A trojan invariant torus for the ER3BP plotted in the $(\xi,\eta)$, $(\dot \xi, \dot \eta)$, $(\xi,\dot \xi)$, and $(\eta, \dot \eta)$ projections. 
    The length $601$ trajectory is given in red, while the computed invariant torus with lines that are constant in coordinates are given in gray.}
    \label{fig:trojan}
\end{figure}

Nearby the stable $\mathrm L_4$ and $\mathrm L_5$ Lagrange points for the unperturbed circular problem, the ER3BP supports invariant tori.
Many asteroids --- called trojans --- are known to orbit along these invariant tori \cite{murray_solar_2000}.
In Fig.~\ref{fig:trojan}, we plot multiple projections of a trajectory initialized on an invariant torus at $\bm \xi = (1/2, \sqrt{3}/2-0.005),$ $\bm p = (0,0)$, nearby the stable $\mathrm L_4$ Lagrange equilibrium point $\bm \xi = (1/2, \sqrt{3}/2)$ for the unperturbed $\epsilon = \mu = 0$ problem.
The plotted trajectory consists of $601$ points, corresponding to a physical time of $7{,}127$ years.
We note that, due to the high-dimensional nature of the problem, different projections of the invariant torus appear significantly different to the eye.
This makes it generally difficult to intuit the high-dimensional geometry of dynamical systems from visuals alone.
However, we do make one observation from Fig.~\ref{fig:trojan}(a) that generalizes well to many other tori: the trojan invariant torus is significantly thinner in some directions than others, a property we generally refer to as ``filamentary'' or ``anisotropic.''
We will find anisotropy is one of the main difficulties of computing invariant tori.

Our goal is to compute an invariant torus $S$ and its rotation vector $\bm \omega$ from a trajectory $\sympmap^t(x)$ for $0\leq t < T$. 
To parameterize the invariant torus, we must define an observable function $\bm h : \mathcal{M} \to \Rbb^D$, $D\geq 2d$, from the state space into Euclidean space so that the torus can be represented as
\begin{equation*}
    \bm h(S(\bm \theta)) = \sum_{\bm k\in \Zbb^d} \bm h_{\bm k} e^{2\pi \I \bm \theta \cdot \bm k}.
\end{equation*}
For the ER3BP, this function is simply the identity $\bm h : (\bm \xi, \bm p) \mapsto (\bm \xi, \bm p)$. 
However, many other smooth maps could be chosen, such as a different coordinate system or a delay embedding.
In the more general case that the underlying state space $\mathcal{M}$ is not Euclidean, the map $\bm h$ allows us to still use Fourier series to represent the geometry.
As long as $\bm h(S)$ is injective, the map $\bm h$ can be pulled back to give an invariant torus on the regular state space.
% If one is only interested in the observable $\bm h$ itself, injectivity is unnecessary. 

When the observable is composed with an invariant torus trajectory, we find that
\begin{equation}
\label{eq:h-invariant-torus}
    \bm h(\sympmap^t(S(0))) = \bm h(S(t \bm \omega)) = \sum_{\bm k} \bm h_{\bm k} \lambda_{\bm k}^t, \qquad \lambda_{\bm k} = e^{2\pi \I \bm \omega \cdot \bm k}.
\end{equation}
The effect is that spatial wavenumbers $\bm k \in \Zbb^d$ of the trajectory are transformed to temporal frequencies $\bm \omega \cdot \bm k \in \Tbb$, where the Diophantine condition guarantees that $\bm \omega \cdot \bm k_1 \neq \bm \omega \cdot \bm k_2$ for $\bm k_1 \neq \bm k_2$.
Given this representation, a more refined statement of our goal is how to recover the Fourier coefficients $\bm h_{\bm k}$ and the rotation vector $\bm \omega$ knowing  $\bm{h}(\sympmap^t(S(0)))$, for $t = 0,1,\dots,N$ for some large $N$.

The process for recovering the torus will proceed in two stages, where first the rotation vector is found, and then the Fourier coefficients are computed. 
However, we will explain the process in the opposite order to motivate the choices we made for computing the rotation vector. 
So, for the rest of this section, assume that $\bm \omega$ is known, and we simply want to find the coefficients $\bm h_{\bm k}$ from the signal.

\subsection{Weighted Birkhoff averages vs.~least squares}
\label{subsec:fourier-coef}
The first of two methods we consider to compute Fourier coefficients is what we refer to as the projection method \cite{blessing_weighted_2024}.
It depends on the inverse Fourier transform formula
\begin{equation}
\label{eq:Fourier-projection}
    \bm h_{\bm k} = \int_{\Tbb^d} \bm h(S(\bm \theta)) e^{- 2\pi\I \bm k\cdot \bm \theta} \dif \bm \theta.
\end{equation}
Because we cannot sample $S$ at will, we connect the quadrature to the data we have, namely a single trajectory.

The connection between quadrature of the form \eqref{eq:Fourier-projection} and trajectories is the Birkhoff ergodic theorem, which states that integrals over invariant sets can be related to time averages over trajectories.
Let $x = S(0)$ be a point on an invariant torus $S$ with a non-resonant rotation vector $\bm \omega$, then the rotation operator $\tau_{\bm \omega} : \bm \theta \mapsto \bm \theta + \bm \omega$ is ergodic on $\Tbb^d$. 
Thus, the Birkhoff ergodic theorem tells us that the finite-time averages of continuous observables $\bm g : \mathcal{M} \to \Rbb^D$ given by
\begin{equation*}
    \mathcal{B}_T[\bm g](x) =  \frac{1}{T}\sum_{t=0}^{T-1} \bm g \circ \sympmap^t(x) = \frac{1}{T} \sum_{t = 0}^{T-1} (\bm g \circ S) \circ \tau_{\bm \omega}^t(0)
\end{equation*}
converge as $T\to\infty$ to a unique limit, defined by the integral
\begin{equation}
    \mathcal B[\bm g](x) = \lim_{T\to \infty}\mathcal{B}_T[\bm g](x) = \int_{\Tbb^d} \bm g(S(\bm \theta)) \dif \bm \theta,
\end{equation}
where we have used the fact that the unique invariant measure of non-resonant rotational dynamics on $\Tbb^d$ is the uniform measure. 

Restricting to the specific case of finding Fourier coefficients, the function we want to average is $\bm g(S(\bm \theta)) = \bm h(S(\bm \theta)) e^{-2\pi\I\bm k \cdot \bm \theta}$. 
Substituting this into the ergodic formula, we can approximate the Fourier coefficients by the average
\begin{equation}
\label{eq:Fourier-time-projection}
    \bm h_{\bm k} \approx \frac{1}{T} \sum_{t=0}^{T-1} \lambda_{\bm k}^{-t}\bm h\circ \sympmap^t(x). 
\end{equation}
As long as we have access to the observation of the trajectory $\bm h \circ \sympmap^t(x)$ and the rotation vector $\bm \omega$ (so that $\lambda_{\bm k}$ is known), we can apply the above ergodic formula to approximate the Fourier coefficients.

The rate of convergence of Birkhoff averages depends on the type of underlying trajectory.
Assuming $x$ is on an invariant torus with a Diophantine rotation vector and the observable and dynamics are sufficiently smooth, it is straightforward to show that $\norm{\mathcal{B}_T[\bm h](x) -\mathcal{B}[\bm h](x)} = \mathcal O(T^{-1})$. 
Fortunately, the rate of convergence for invariant tori can be accelerated by appropriately weighting the sum as
\begin{equation*}
% \label{eq:WBA}
    \mathcal{WB}_T[\bm g](x) = \sum_{t = 0}^{T-1} w_t \bm g \circ \sympmap^t (x),    
\end{equation*}
where the weights $w_t$ are sampled from a compactly supported function $w\in C^{\infty}([0,1])$ function via
\begin{equation*}
    w_t = \left[\sum_{s=0}^{T-1} w\left(\frac{s+1}{T+1} \right) \right] w\left(\frac{t+1}{T+1}\right), \qquad w(x) = \exp\left(-\frac{1}{x(1-x)}\right).
\end{equation*}
This technique is known as a weighted Birkhoff average. 
The error $\norm{\mathcal{WB}_T[\bm g](x)-\mathcal{B}[\bm g](x)}$ for $\bm h,$ $S,$ and $\mathcal{M}$ in $C^M$ converges as $\mathcal O(T^{-m})$ provided that $M > d + m\nu$ \cite{Das2017}, and the error for the analytic system can converge as $\mathcal O(e^{- r T^\xi})$ for some rates $r>0$ and $0 < \xi \leq 1$, depending on the dimension and the torus \cite{Tong2024}.

We note that ergodic averages of chaotic trajectories converge at a slower convergence rate than for invariant tori \cite{sander2020}. 
The empirical rate is approximately $\mathcal O(T^{-1/2})$ regardless of the weighting, where the exponent is due to the heuristic that a chaotic average should converge approximately according to a central limit theorem.
The discrepancy of rates between chaos and tori allows for the weighted Birkhoff average to efficiently classify trajectories.
While sorting trajectories is not the focus of this paper, this is a major motivation for computing invariant tori from trajectories.

We note that the validity of the projection formula \eqref{eq:Fourier-time-projection} relies upon an assumption of orthogonality between Fourier modes.
This is a good assumption in the context of spatial averages, as 
\begin{equation}
\label{eq:Fourier-orthogonality}
    \ip{e^{2\pi \I \bm k \cdot \star}}{e^{2\pi \I \bm k' \cdot \star}} = \int_{\Tbb^d} e^{2\pi \I (\bm k' - \bm k)\cdot \bm \theta} \dif \bm \theta = \begin{cases}
        1 & \bm k = \bm k', \\
        0 & \bm k \neq \bm k'.
    \end{cases}
\end{equation}
This orthogonality relationship does not exactly hold for the finite-time ergodic average, however, where the approximate integral
\begin{equation*}
    \ip{e^{2\pi \I \bm k \cdot \star}}{e^{2\pi \I \bm k' \cdot \star}} \approx \sum_{t=0}^{T-1} w_t \lambda_{\bm k}^{-t} \lambda_{\bm k'}^t
\end{equation*}
is nonzero for all $\bm k \neq \bm k'$.
Moreover, this error increases as the difference between frequencies grows, owing to the fact that higher Fourier coefficients are increasingly non-smooth.

Thinking of the weighted Birkhoff average as a quadrature rule on $\Tbb^d$, it is informative to compare to the trapezoidal rule.
The trapezoidal rule is a Gaussian quadrature rule $\Tbb$, meaning its tensor product rule exactly preserves the orthogonality relation \eqref{eq:Fourier-orthogonality} for $\norm{\bm k - \bm k'}_\infty < 2K_{\mathrm{trap}}-1$, where $K_{\mathrm{trap}}$ grid points in each direction \cite{trefethen_exponentially_2014}.
The projection that one obtains from the trapezoidal rule is the orthogonal discrete Fourier transform, which is what allows the continuous projection rule to also be a discrete one.
On the other hand, because the weighted Birkhoff average is not orthogonal, the continuous projection formula will only be an approximation to the true discrete projection

To tackle the non-orthogonality, we use the least-squares approach
\begin{equation}
\label{eq:Fourier-lsqr}
    \min_{\{\bm h_{\bm k_n}\}} \sum_{t=0}^{T-1}c_t \norm{\left(\sum_{n=1}^N \bm h_{\bm k_n}\lambda_{\bm k_n}^t\right) - \bm h(\sympmap^t(x)) }^2 = \min_{X} \norm{A X - B}_W^2,
\end{equation}
where the weights $c_t$ are either ergodic average $c_t = 1/T$ or the weighted average $c_t = w_t$, $A\in \Cbb^{T\times N}$ is the matrix representing the $N$ Fourier modes $\bm k_n \in \Zbb^d$
\begin{equation}
\label{eq:A-Rh}
    A = \begin{pmatrix}
        | & & | \\
        \bm \Lambda_{\bm k_1} & \dots & \bm \Lambda_{\bm k_N} \\
        | & & |
    \end{pmatrix}, \qquad
    \bm \Lambda_{\bm k} = \begin{pmatrix}
        \lambda_{\bm k}^0 \\ \vdots \\ \lambda_{\bm k}^{T-1}
    \end{pmatrix},
\end{equation}
the matrices $X\in\Cbb^{N\times D}$ and $B\in \Rbb^{T \times D}$ are the frequency and time representations of the signal
\begin{equation*}
    X = \begin{pmatrix}
        \bm h_{\bm k_1}^\trans \\ \bm h_{\bm k_2}^\trans \\ \vdots \\ \bm h_{\bm k_N}^\trans
    \end{pmatrix}, \qquad B = \begin{pmatrix}
        \bm h(\sympmap^0(x))^\trans \\ \bm h(\sympmap^1(x))^\trans \\ \vdots \\ \bm h(\sympmap^{T-1}(x))^\trans
    \end{pmatrix},
\end{equation*}
and $\norm{\star}_W$ is the weighted norm
\begin{equation*}
    \norm{M}^2_W = \tr(M^* W M), \qquad W = \Diag(c_0, \dots, c_{T-1}),
\end{equation*}
where $M^*$ is the conjugate transpose of $M$. 
The least-squares system solution can be represented exactly using the normal matrix formula as
\begin{equation*}
    X = (A^* W A)^{-1} (A^* W B).
\end{equation*}
In the formula above, the (weighted) Birkhoff average for the modes $\bm k_n$ is simply the $A^* W B$ term.
Then, the least-squares solution corrects the non-orthogonality of that approach by inverting the normal matrix $A^* W A$. 
(In practice, we solve the system by performing a QR factorization $QR = W^{1/2} A$ and computing $X = R^{-1} Q^* W^{1/2} B$.)

Now that we have two methods --- projection and least-squares --- we return to the question of how many modes $\bm k_n$ to compute from the signal. 
In 1D, a partial answer to this question was given in \cite{ruth_finding_2024}, where the condition number of the least-squares problem was used to determine the size. 
However, the anisotropy of the higher-dimensional problem requires resolutions in different directions for the best interpolation, which the condition number argument is insufficient for.
This issue compounds the already mentioned issue of a lack of orthogonality, meaning increasing resolution can worsen results.

To demonstrate the issue, we return to our example trajectory in the ER3BP.
Using the method of Section \ref{sec:rotation-vector}, we find the rotation vector to be $\bm \omega = (0.00344, 0.07921)$, where the first and second entries correspond to rotating the ``short'' and ``long'' ways around the torus respectively.
Using this frequency, we compute Fourier coefficients from the projection and least-squares method, using both the unweighted and the weighted versions.
To evaluate the residual, the methods are all trained on the first $T=550$ points of the trajectory for wavenumbers $\abs{k_j} \leq K_j$, where we scan different resolutions $0\leq K_1 \leq 8$ and $0\leq K_2 \leq 25$.
For each method and resolution $(K_1,K_2)$, we form an approximate invariant torus of the form
\begin{equation*}
    \hat{\bm h}(\bm \theta) \approx \bm h \circ S(\bm \theta) = \sum_{k_1=-K_1}^{K_1} \sum_{k_2=-K_2}^{K_2} \bm h_{\bm k} e^{2\pi \I \bm k \cdot \bm \theta}.
\end{equation*} 
This torus is validated by computing the relative error between the predicted dynamics on the next $M = 51$ points on the trajectory and the real dynamics as
\begin{equation}
\label{eq:validation-error}
    R_h^2 = \frac{1}{R_{h0}^2} \sum_{t=T+1}^{T+M} \norm{\bm h(\sympmap^t(x)) - \hat{\bm h}(\bm \theta + \bm \omega t)}^2, \qquad R_{h0}^2 = \sum_{t=T+1}^{T+M}\norm{\bm h(\sympmap^t(x))}^2.
\end{equation}

\begin{figure}
    \centering
    \includegraphics[width=0.8\linewidth]{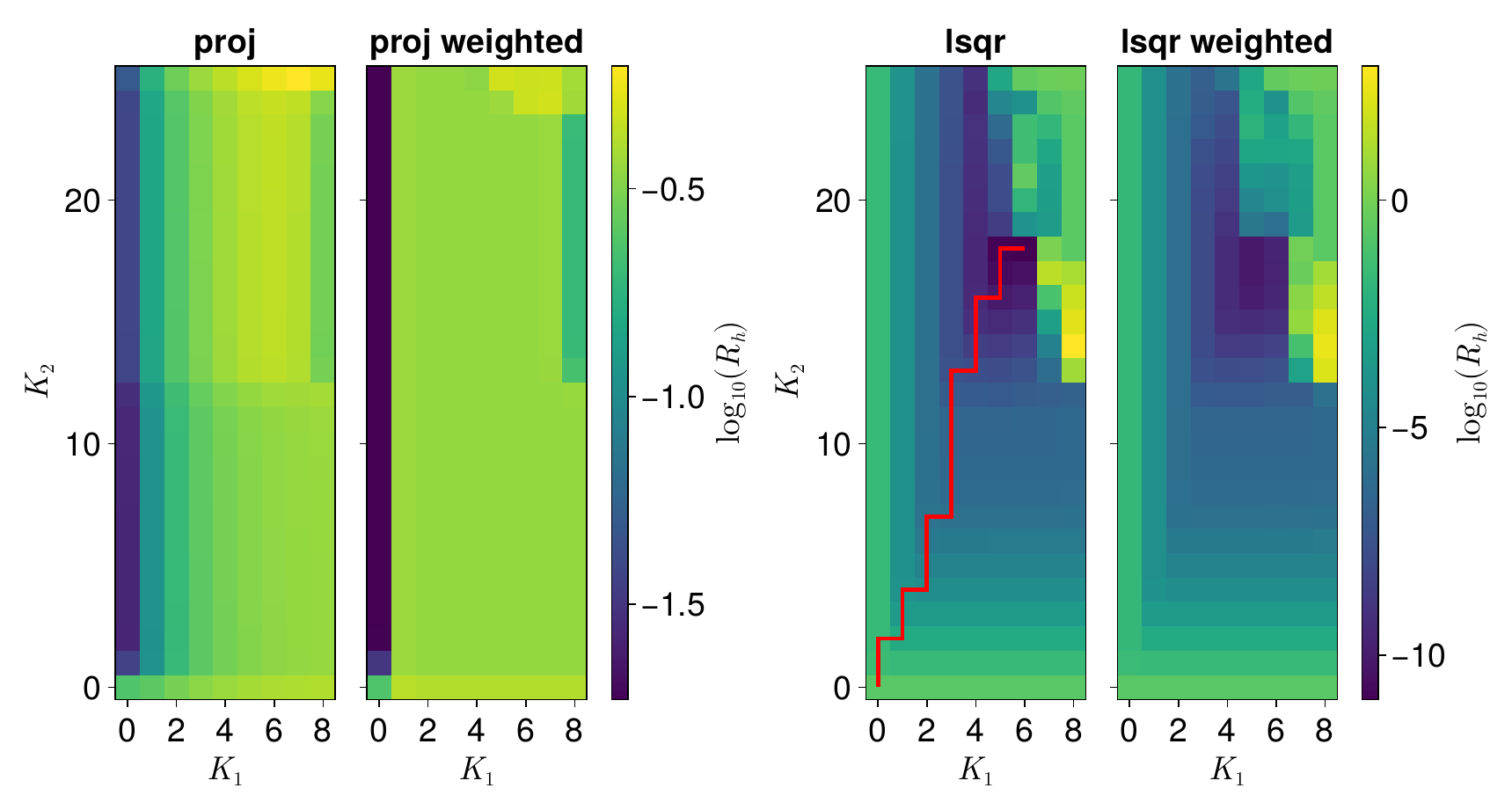}
    \caption{The validation error \eqref{eq:validation-error} for the (left) projection and (right) least-squares methods, for both the unweighted and weighted versions of each method. The red path on the unweighted least-squares method indicates the path of the adaptive method (see Sec.~\ref{subsec:fourier-adaptive}).}
    \label{fig:val}
\end{figure}

\begin{table}
    \centering
    \begin{tabular}{c|c|c|c|c}
         Method                 & $L$                                            & Best Error $R_h$        & $K_1$ & $K_2$ \\ \hline \hline
         Projection             & $I$                                            & $2.67\times 10^{-2}$  & $0$   & $2$   \\
         Weighted Projection    & $I$                                            & $1.87\times 10^{-2}$  & $0$   & $24$  \\
         Least Squares          & $I$                                            & $1.07\times 10^{-11}$ & $6$   & $18$  \\
         Weighted Least Squares & $I$                                            & $7.86\times 10^{-11}$  & $5$   & $18$  \\ 
         Least Squares          & $\begin{pmatrix} 1 & 0 \\ 1 & 1 \end{pmatrix}$ & $1.61\times 10^{-6}$  & $8$  & $8$  \\
         Least Squares          & $\begin{pmatrix} 1 & 0 \\ 2 & 1 \end{pmatrix}$ & $9.23\times 10^{-6}$  & $14$  & $5$   \\
         Least Squares          & $\begin{pmatrix} 2 & 1 \\ 1 & 1 \end{pmatrix}$ & $3.64\times 10^{-5}$  & $5$   & $10$  \\
    \end{tabular}
    \caption{Best validation errors and their resolutions. The first four lines match Fig.~\ref{fig:val}. The final three lines give the unweighted least-squares validation error for three valid rotation vector transformations $\bm \omega' = L\bm \omega$.}
    \label{tab:val}
\end{table}

In Fig.~\ref{fig:val}, we plot the validation error as a function of $K_1$ and $K_2$ for each coefficient recovery method.
Then, in the first four rows of Table~\ref{tab:val}, we report the optimal validation error, along with the optimal values of $K_1$ and $K_2$.
We find that the least-squares methods for finding Fourier coefficients outperform the projection methods by nine orders of magnitude, with the best performance being the $(K_1,K_2) = (6,18)$ unweighted least squares. 

There are two main takeaways from this numerical experiment.
The first takeaway is that for the trajectory length of $550$ iterations, the Fourier modes for both the projection and weighted projection methods are far from being orthogonal, leading to the drastic errors. 
This can be attributed to the small value of $\omega_1 = 0.00344$, as there are only $3$ full periods of $e^{2\pi \I \omega_1 t}$ over the course of the considered trajectory. 
As the trajectory length increases, the weighted projection would converge much faster than the unweighted projection method, but this is only a limiting result, and cannot necessarily be applied for ``short'' trajectories.
One could fix this problem by running a longer trajectory, but the difference in integration time can easily exceed the difference in computation time between the projection and least-squares methods.
Moreover, this trajectory length is sufficient to find the rotation vector to high accuracy, so we do not need more information for that part of the problem.

The second takeaway is that the least-squares method can achieve very accurate results regardless of the weighting, but only if the anisotropy is accounted for.
Clearly, if we had only considered projections where $K_1 = K_2$, a validation error on the order of $10^{-10}$ could never have been achieved, as the high Fourier modes in the $\theta_1$ direction would overfit while the modes in the $\theta_2$ direction would underfit. 

\subsection{Adaptivity in the number of Fourier modes}
\label{subsec:fourier-adaptive}
While we have shown that an anisotropic Fourier representation of the invariant torus can achieve good validation error, we have not yet described how to determine the correct resolution.
To do so, consider a length $T$ trajectory of observations $\bm h(\sympmap^t(x))$, the rotation vector $\bm \omega$, and a validation fraction $0 < \gamma < 1$.
% To address this, we consider a greedy approach to computing the best resolution.
To address this, we attempt to solve the following problem:
\begin{equation}
\label{eq:torus-residual}
    \min_{\bm K} R_{h}^2[\hat{\bm h}], \qquad \hat{\bm h}(\bm \theta) = \sum_{-\bm K \preceq \bm k \preceq \bm K} \bm h_{\bm k} e^{2\pi \I \bm k \cdot \bm \theta},
\end{equation}
where $\bm K = (K_1, \dots, K_d)$ is the Fourier resolution of the approximate torus $\hat{\bm h}$, determined by solving the unweighted least-squares system \eqref{eq:Fourier-lsqr} on the first $\floor{(1-\gamma)T}$ elements of the trajectory, where $\floor{\star} : \Rbb\to \Zbb$ is the floor function.
The objective for each value is computed by Eqs.~\eqref{eq:validation-error} and \eqref{eq:torus-residual}, computed on the final $1-\floor{(1-\gamma)T}$ points on the trajectory.
In this paper, we fix $\gamma = 0.05$.
In essence, the minimization solves for the invariant torus that best predicts the dynamics. 

The method of optimization is simple: we perform a discrete descent on the validation error (see Fig.~\ref{fig:val}).
For the optimization, there are two hyperparameters: a ``hill tolerance'' $\eta\geq 1$ and a maximum number of Fourier coefficients $K_{\mathrm{max}}\leq \floor{(1-\gamma)T}$.
The descent is initialized with a constant interpolant, corresponding to $K_j = 0$ for all $1\leq j \leq d$.
To determine the next resolution $\bm K$, we increase the resolution in each coordinate and compute the residual $\tilde{R}_h^{(j)} = \tilde{R}_h^{(j)}(\bm K + \bm e_j)$ for $1 \leq j \leq d$, where $\bm e_j \in \Zbb^d$ is the unit vector in the $j$ direction.
After computing the residuals, we move in the direction of steepest descent $j = \argmin_{j'} \tilde{R}_h^{(j')}$ (note that residual need not decrease). 
Along the descent, we keep track of the value of $\bm K = \bm K^*$ that gives the smallest residual $R_h(\bm K^*)$.
We stop the descent when there is no candidate direction such that the residual improves $\tilde{R}_h(\bm K + \bm e_j) < \eta \tilde{R}_h(\bm K^*)$ and the resolution is bounded by $K=\prod_j (2 K_j+1) \leq K_{\mathrm{max}}$.

In the case that $\eta=1$, the descent stops when any increase in resolution would increase the residual.
However, we have found that this can stop the descent in local minimum, while the global minimum may have a higher resolution.
In these cases, values of $\eta > 1$ allow the descent algorithm to traverse ridges in the validation error where it could possibly find better minima.
There is no guarantee that this algorithm finds the global minimum of the validation error, but choosing $\eta > 1$ (typically we choose $\eta = 10$) allows for a better chance.

The primary benefit of this algorithm is that it is practical.
A worst-case asymptotic cost of this algorithm is $\Ocal(d T K_{\mathrm{max}}^2)$, where $d$ is the dimension of the torus and $T$ is the length of the trajectory.
This cost is dominated by the computation of $R_h(\bm K + \bm e_j)$ at each step, where the solutions of the least-squares problems are computed by QR update formulas at each step (see Appendix \ref{app:qr-updates} for details). 
However, the asymptotic bound is pessimistic because the stopping condition is often met before the worst-case time.
More often, if $\bm K$ is the resolution of the largest least-squares problem solved and $K = \prod_j (2K_j+1) \leq (1-\gamma) T$ is the width of the least-squares problem, then the cost is given by $\Ocal(d T K^2)$. 

Better results could be obtained, e.g., by performing a search over all possible values of $\bm K$, but this would be significantly more expensive.
This may be practical if $\sympmap$ is very expensive to compute.
Otherwise, the value of \eqref{eq:torus-residual} can be reduced much more simply by computing a longer trajectory. 
Moreover, the parameterization method \cite{haro2016} can be used to iteratively improve an initial guess of an invariant torus when the found torus needs to be refined (this also gives an avenue to prove the existence of the invariant torus, see \cite{Figueras2017}).

\subsection{Non-uniqueness of the rotation vector}
\label{subsec:nonuniqueness}
We have seen from the ER3BP example that the least-squares method of computing Fourier coefficients can dramatically outperform the weighted Birkhoff average projection method.
However, up to now, there is an issue that we have ignored: the parameterization of an invariant torus is not unique.

To investigate the non-uniqueness, consider two parameterizations of the same invariant torus $S$ and $S'$ with non-resonant rotation vectors $\bm \omega$ and $\bm \omega'$. 
Then, it must be the case that $S' = S\circ \bm G$ where $\bm G : \Tbb^d \to \Tbb^d$ is a automorphism of the torus regarded as a smooth manifold.
Automorphisms of the torus can be decomposed as
\begin{equation*}
    \bm G = L \circ \bm G_0,
\end{equation*}
where $\bm G_0$ induces the identity on the singular homology $H_1(\Tbb^d) \approx \Zbb^d$ and $L \in GL(d,\Zbb)$ induces a non-trivial isomorphism on $H_1(\Tbb^d)$ when $L\neq I$.
Here, $GL(d,\Zbb)$ is the group of integer matrices $L \in \Zbb^{d\times d}$ with integer-valued inverses $L^{-1} \in \Zbb^{d\times d}$. The inverse condition is satisfied if and only if $\det L = \pm 1$.
By requiring that $S \circ \tau_{\bm \omega} = S\circ \bm G \circ \tau_{\bm \omega'}$, one finds that $\bm \omega = L \bm \omega'$ and $\bm G_0(\bm \theta) = \bm \theta + \bm \Delta$ for some constant $\bm \Delta\in\Tbb^d$.
This tells us that the non-uniqueness in torus representation can be expressed by two quantities: $\bm \Delta$ for the choice of origin and the $L$ for the choice of homology generators.

The shift $\bm \Delta$ is easily fixed by requiring $\bm h(S(0)) = \bm h(x)$, but the best choice of $L \in GL(d,\Zbb)$ is less intuitive.
To see the numerical implications of the choice of $L$, we first recognize that Fourier coefficients change under the automorphism as
\begin{align*}
    \bm h(S'(\bm \theta')) &= \sum_{\bm k'} \bm h_{L^\trans \bm k'} e^{-2\pi \I \bm k' \cdot \bm \theta'},
\end{align*}
where $\bm k' = L^{-\trans} \bm k$.
In this way, the choice of $L$, and therefore the choice of $\bm \omega'$, can be identified with a choice of the wavenumber ordering.

\begin{figure}
    \centering
    \includegraphics[width=0.8\linewidth]{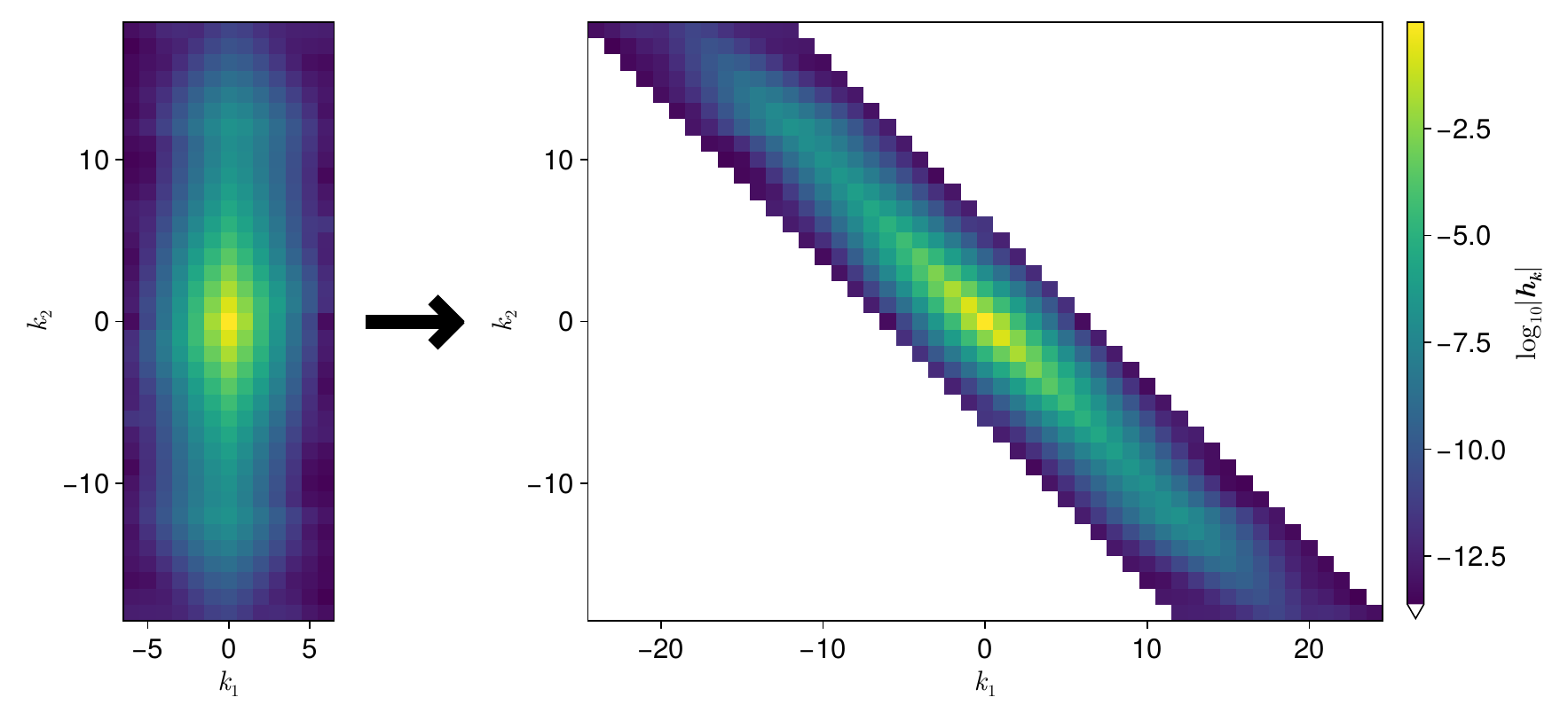}
    \caption{A schematic of how the Fourier coefficients for the ER3BP example shear under the transformation $L_1$.}
    \label{fig:Fourier-shear}
\end{figure}
Within this non-uniqueness lies a problem: a Fourier representation for one parameterization of an invariant torus might lie neatly within a rectangular domain of wavenumbers, while the transformed torus may have a skewed representation. 
To see this, recall the ER3BP example and let 
\begin{equation*}
    L_1 = \begin{pmatrix}
        1 & 0 \\ 1 & 1
    \end{pmatrix}.
\end{equation*}
(We immediately recognize that $L_1\in GL(2,\Zbb)$, as its inverse is in $\Zbb^{2\times 2}$ and it has determinant $1$.)
In Fig.~\ref{fig:Fourier-shear}, we show how the Fourier coefficients from the original parameterization change under this transformation.
We see that what was once a vertical exponentially decaying function is transformed to be diagonal.
The upshot is that a resolution of $(K_1,K_2) = (24,18)$ would be needed to capture all of the same Fourier modes for for the sheared torus that are captured in a resolution of $(K_1,K_2) = (6,18)$ for the original.
So, as long as we compute the Fourier coefficients over a rectangle, it is preferable to choose the homology generators with the most compact representation.

To make this observation precise, in Table \ref{tab:val}, we scan over $K_1$ and $K_2$ to find the best validation error \eqref{eq:validation-error} for the transformed torus.
We perform this scan for the transformation $L_1$, as well as two other matrices
\begin{equation*}
    L_2 = \begin{pmatrix}
        1 & 0 \\ 2 & 1
    \end{pmatrix}, \qquad L_3 = \begin{pmatrix}
        2 & 1 \\ 1 & 1
    \end{pmatrix}.
\end{equation*}
In all three cases, the best validation error reduces by $4$ orders of magnitude, showing that the choice of homology can significantly change the quality of the computed invariant torus.
In Section \ref{subsec:refine-rotation}, we will return to this detail on parameterization to determine a ``best'' rotation vector. 
\section{Computing the rotation vector}
\label{sec:rotation-vector}
Assuming the rotation vector is known, we have seen that it is possible to obtain an accurate Fourier representation of an invariant torus from a trajectory.
Now, we turn our attention to finding the rotation vector itself.
Before explaining the details in Secs.~\ref{subsec:BirkhoffRRE}, \ref{subsec:initial-rotation}, and \ref{subsec:refine-rotation}, we show an example of why it is more difficult to compute high-dimensional rotation vectors than in one dimension.

\begin{figure}
    \centering
    \includegraphics[width=0.9\linewidth]{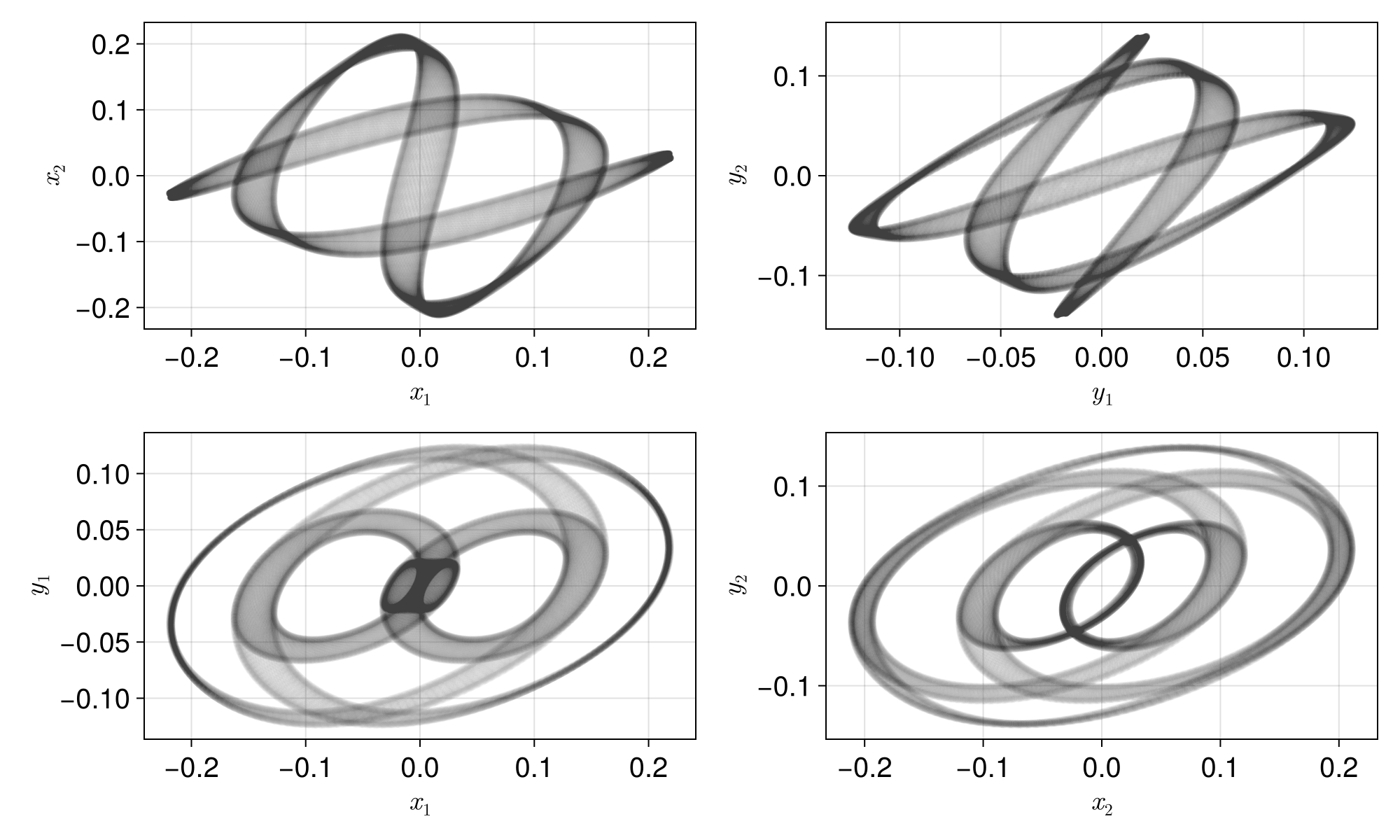}
    \caption{Four projections of a trajectory of the coupled standard map \eqref{eq:coupled-sm}. }
    \label{fig:wonky-standard-map}
\end{figure}

For the example, consider the coupled standard map
\begin{equation}
\label{eq:coupled-sm}
    \bm x_{t+1} = \bm x_t + \bm y_{t+1}, \qquad \bm y_{t+1} = \bm y_t + \frac{1}{2\pi} K_{\mathrm{sm}} \sin(2\pi \bm x_t),\qquad K_{\mathrm{sm}} = \begin{pmatrix}
        0.4 & 0.2 \\ 0.2 & 0.5
    \end{pmatrix},
\end{equation}
where $\bm x \in \Tbb^2$, $\bm y \in \Rbb^2$ and the function $\sin(\cdot)$ is taken elementwise.
The symplectic map on the state space $\Tbb^2 \times \Rbb^2$ is defined as $\sympmap: (\bm x_t, \bm y_t)\mapsto(\bm x_{t+1}, \bm y_{t+1})$.
In Fig.~\ref{fig:wonky-standard-map}, we plot four projections of a trajectory starting at $\bm x = (0.1,0.1)$ and $\bm y = (0.1,0.01)$.
The rotation vector is $\bm \omega = (0.02487, 0.02414)$, as is determined through the methods within this section.

There are three things that make this torus specifically difficult:

\paragraph{1. The torus is far from perturbative.}
If $(x_1,y_1)$ was uncoupled from $(x_2,y_2)$ (i.e.~$K_{\mathrm{sm}}$ is diagonal), the 2D invariant torus would be a tensor product of two 1D invariant tori.
So, when plotted in the $(x_1,y_1)$ and $(x_2,y_2)$ planes, we would observe a non-self-intersecting projection.
This is not the case we see in Fig.~\ref{fig:wonky-standard-map}; it is clear that the dynamics of $(x_1,y_1)$ is strongly coupled with the dynamics of $(x_2,y_2)$.
This makes the question of continuation to this torus very difficult.
It may be possible to continue from a tensor product of 1D tori or even the $K_{\mathrm{sm}}=0$ system, but this would require a method that is robust to resonant bifurcations.
    
\paragraph{2. Winding arguments break down.}
Consider an invariant torus from the coupled standard map \eqref{eq:coupled-sm} where $\bm y \in \Rbb^2$ is a graph over $\bm x \in \Tbb^2$.
In this case, one can take averages of the ``momentum'' observable $ \bm h_y : (\bm x,\bm y) \mapsto \bm y$ to find the rotation vector:
\begin{equation*}
    \bm \omega = \mathcal B[\bm h_y](\bm x_0, \bm y_0).
\end{equation*}
However, the example torus in Fig.~\ref{fig:wonky-standard-map} does not have the graph property.
In some cases this can be overcome by a simple coordinate transformation.
For this example, the complicated high-dimensional geometry of the torus in Fig.~\ref{fig:wonky-standard-map} makes it difficult to intuit or compute such a coordinate transformation. 

\paragraph{3. The torus is dominated by high-wavenumber Fourier modes. }
One seemingly straightforward method to compute rotation numbers is by a discrete Fourier transform (DFT). 
If the trajectory's signal is truly of the form \eqref{eq:h-invariant-torus}, then one expects that a DFT will have peaks near the frequencies $\Omega_{\bm k} = \bm \omega \cdot \bm k$.
So, in principle, the peaks of that Fourier transform have all of the necessary information.

\begin{figure}
    \centering
    \includegraphics[width=1.0\linewidth]{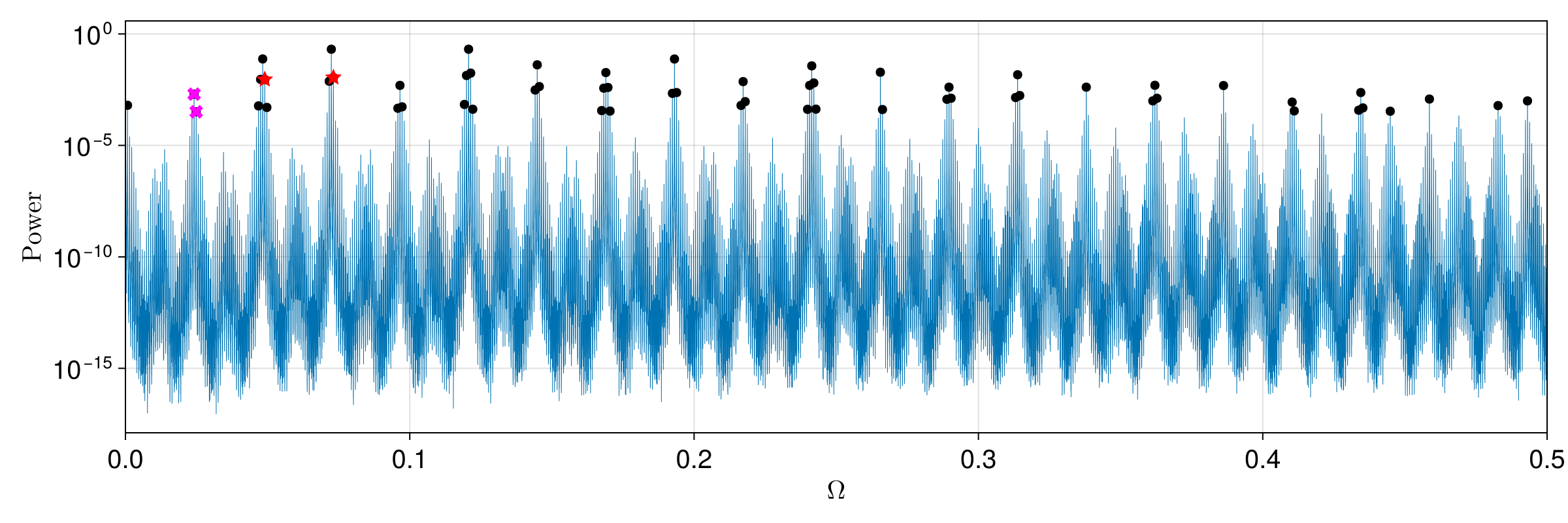}
    \caption{A $C^\infty$-windowed discrete Fourier transform of a length $10000$ trajectory of the 2D standard map \eqref{eq:coupled-sm} shown in Fig.~\ref{fig:wonky-standard-map} with the observable function in Eq.~\eqref{eq:sm-observable}. 
    The 60 highest peaks $(\Omega_j, H_j)$ for $1\leq j \leq 60$ are plotted in black dots, cf.~Table \ref{tab:wonky-coefficients}.
    The only two peaks in the 15 largest that form a valid rotation vector are marked with red stars ($j=12,14$).
    The optimal rotation vector for the anisotropy of the invariant torus is marked in magenta `x's ($j=31,60$).}
    \label{fig:wonky-spectrum}
\end{figure}

For a first attempt at a DFT, let $\bm h : \Tbb^2 \times \Rbb^2 \to \Rbb^4$ be an observable function defined by
\begin{equation}
\label{eq:sm-observable}
    \bm h(\bm x, \bm y) = \begin{pmatrix}
        (y_1 + R) \cos(2\pi x_1) \\
        (y_1 + R) \sin(2\pi x_1) \\
        (y_2 + R) \cos(2\pi x_2) \\
        (y_2 + R) \sin(2\pi x_2)
    \end{pmatrix},
\end{equation}
where $R=0.5$.
Then, we compute a $C^\infty$-windowed fast Fourier transform of a length $10000$ trajectory of observables.
The results of the DFT are shown in Fig.~\ref{fig:wonky-spectrum}, where the spectrum appears to be structured with many peaks.
However, this still leaves two questions:
\begin{enumerate}
    \item[(a)] Are the DFT frequencies accurate enough? 
    \item[(b)] Which peaks correspond to a valid rotation vector?
\end{enumerate}

For Question (a), the answer is no. 
A length $T$ trajectory yields a DFT frequency resolution of $\Ocal(T^{-1})$, which is not enough to practically achieve the parameterization accuracies of Sec.~\ref{sec:parameterization}. 
In Sec.~\ref{subsec:BirkhoffRRE}, we solve this problem using an extrapolation method --- Birkhoff reduced rank extrapolation (Birkhoff RRE) --- which has the benefit of computing many frequencies $\Omega_j \approx \bm \omega \cdot \bm k_j$ simultaneously to high accuracies.
We then project the frequency modes back onto the trajectory, giving the coefficients $\bm H_j \approx \bm h_{\bm k_j}$ ordered by the magnitudes $H_1 \geq H_2 \geq \dots \geq H_{J^*}$ where $H_j = \norm{\bm H_j}$.
In this way, the invariant torus is represented as
\begin{equation*}
    \bm h(\sympmap^t(\bm x)) = \sum_{j=-J^*}^{J^*} \bm H_j e^{2\pi \I \Omega_j t},
\end{equation*}
where we use the convention that $\Omega_{-j} = -\Omega_j$ and $\Omega_0 = 0$.

For Question (b), we note that for 1D invariant tori in 2D Euclidean space, the answer is almost always simple: choose the highest peak in the spectrum to be the rotation number (i.e.~$\omega = \Omega_1$). 
This works because 1D tori embedded in 2D space cannot self-intersect, meaning that the $k=1$ Fourier mode is almost always the most prominent.
Because 2D tori in 4D space are not co-dimension 1, they are free to wind around themselves in complicated ways, meaning this reasoning breaks down for higher dimensions.
It is still the case that our desired rotation vector corresponds to two peaks, marked by the magenta `x's in Fig.~\ref{fig:wonky-spectrum}, but it is not obvious why one would choose those peaks from the power spectrum alone.
The best frequencies are far from the most prominent, and essentially resemble the rest of the observed frequencies, a property that is typical of filamentary invariant tori.

\begin{table}
    \centering
    \begin{tabular}{c|c|c|c|c}
        $j$ & $H_j = \norm{\bm H_j} = \norm{\bm h_{\bm k_j}}$ & $\Omega_j = \bm \omega\cdot \bm k_j$ & $(\bm k_j)_{1}$ & $(\bm k_j)_{2}$ \\ \hline
        1 & 2.0577-01 & 0.120688 & 0 & 5 \\
        2 & 2.055901e-01 & 0.072413 & 0 & 3 \\
        3 & 7.535484e-02 & 0.048275 & 0 & 2 \\
        4 & 7.485108e-02 & 0.193101 & 0 & 8 \\
        5 & 4.106522e-02 & 0.144826 & 0 & 6 \\
        6 & 3.691528e-02 & 0.241376 & 0 & 10 \\
        7 & 1.913862e-02 & 0.265514 & 0 & 11 \\
        8 & 1.845954e-02 & 0.168963 & 0 & 7 \\
        9 & 1.747752e-02 & 0.121424 & 1 & 4 \\
        10 & 1.475290e-02 & 0.313789 & 0 & 13 \\
        11 & 1.368350e-02 & 0.119953 & 1 & -6 \\
        \textbf{12} & \textbf{1.113079e-02} & \textbf{0.073148} & \textbf{1} & \textbf{2} \\
        13 & 9.119838e-03 & 0.047540 & 1 & -3 \\
        \textbf{14} & \textbf{9.043146e-03} & \textbf{0.049011} & \textbf{1} & \textbf{1} \\
        15 & 7.556361e-03 & 0.071677 & 1 & -4 \\
        \vdots & \vdots & \vdots & \vdots\\
        31 & 1.986039e-03 & 0.024138 & 0 & 1 \\
        \vdots & \vdots & \vdots & \vdots \\
        60 & 3.211915e-04 & 0.024873 & 1 & 0
    \end{tabular}
    \caption{Locations and magnitudes of the peaks of the Fourier spectrum in Fig.~\ref{fig:wonky-spectrum}.
    In the final two columns, the corresponding wavenumbers of the modes are reported.
    The redundant complex conjugate peaks are not reported (i.e.~we only report $k_1 \geq 0$, and $k_2 > 0$ if $k_1=0$)}
    \label{tab:wonky-coefficients}
\end{table}

In Sec.~\ref{subsec:initial-rotation}, we start by solving the simpler problem of finding a valid rotation vector $(\Omega_{j_1}, \Omega_{j_2})$ from the measured frequencies $\Omega_j$.
However, it is still a nontrivial question as to when two frequencies are a valid rotation vector.
To see this, recall from Sec.~\ref{subsec:nonuniqueness} that any valid rotation vector $(\Omega_{j_1}, \Omega_{j_2})$ can be represented by a transformation $L\in GL(d,\Zbb)$ to the optimal rotation vector as
\begin{equation*}
    \begin{pmatrix}
        \Omega_{j_1} \\ \Omega_{j_2}
    \end{pmatrix} = L \bm \omega,
\end{equation*}
where $\det L = \pm 1$. 
We know the observed frequencies $\Omega_j = \bm \omega \cdot \bm k_j$ for some wavenumbers $\bm k_j$, so we can identify the matrix $L$ as
\begin{equation*}
    L = \begin{pmatrix}
        \bm k_{j_1}^\trans \\ \bm k_{j_2}^\trans
    \end{pmatrix}.
\end{equation*}
By checking whether $\det L = \pm 1$ (and therefore it is in $GL(d,\Zbb)$), this gives us a simple \textit{a posteriori} test to determine whether a rotation vector is valid.

In Table \ref{tab:wonky-coefficients}, we report the powers, frequencies, and the (so far undiscovered) wavenumbers of the 15 highest peaks in Fig.~\ref{fig:wonky-spectrum}, computed by the methods of Sec.~\ref{subsec:BirkhoffRRE}.
We notice that the 8 most prominent frequencies $\Omega_j$ for $1\leq j\leq 8$ are all unusable as elements of a valid rotation vector, as $L$ cannot have a determinant of $\pm 1$ when $\bm k_{j_\ell} = (0,n)$ for $n\geq 2$. 
In fact, within the 15 highest peaks, there is only a single choice of $(j_1,j_2) = (12,14)$:
\begin{equation*}
    \begin{pmatrix}
        \Omega_{12} \\ \Omega_{14}
    \end{pmatrix} = \begin{pmatrix}
        \bm k_{12}^\trans \\ \bm k_{14}^\trans
    \end{pmatrix}\bm \omega = \begin{pmatrix}
        1 & 2 \\ 1 & 1
    \end{pmatrix} \bm \omega.
\end{equation*}
We can confirm the above matrix has determinant $-1$, meaning $(\Omega_{12},\Omega_{14})$ is a valid rotation vector.

In Sec.~\ref{subsec:initial-rotation}, we use a Bayesian maximum \textit{a posteriori} (MAP) estimation problem to determine a candidate rotation number and the associated wavenumbers.
As input, we use the output data from Birkhoff RRE of $\bm \Omega = (\Omega_1, \dots, \Omega_{\JMAP})$ and powers $\bm H = (H_1, \dots, H_{\JMAP})$ for ${\JMAP} < J^*$.
Then, the MAP estimation problem is
\begin{equation}
\label{eq:map}
    \max_{\bm \omega, \bm \pi} \prob(\bm \omega, \bm \pi \mid \bm \Omega, \bm H ),
\end{equation}
where $\bm \pi : j \mapsto \bm k_j$ is a candidate wavenumber labeling. 
The prior of the inference problem will be informed by the assumed smoothness of the invariant torus, which encodes the belief that lower wavenumbers $\norm{\bm k_j}$ should be associated with larger coefficients $H_j$. 

The third step in Sec.~\ref{subsec:refine-rotation} is to choose an optimal rotation vector for parameterization. 
The rough idea of the algorithm is to compute a transformation $L\in GL(d,\Zbb)$ where $\bm \omega' = L \bm \omega$ has the property that $\omega'_1$ points in a shortest direction around the torus, $\omega'_2$ in the second shortest direction, etc.
This notion is made precise by Korkine-Zolatarev reduced lattice basis, which can be quickly computed for small $d$. 

Finally, in Sec.~\ref{subsec:islands}, we describe how the overall algorithm changes to account for trajectories on island chains.

\subsection{Overview of Birkhoff RRE}
\label{subsec:BirkhoffRRE}
The primary purpose of Birkhoff RRE is to accurately translate the trajectory into frequency information.
This translation occurs by way of finding $2J+1$ optimal coefficients $\bm c = (c_{-J}, \dots, c_J) \in \Rbb^{2J+1}$ for the Birkhoff RRE average
\begin{equation}
\label{eq:WBA}
    \mathcal{BR}[\bm h](\bm x) := \sum_{j=-J}^J c_j \bm h(\sympmap^{j+J}(\bm x)).
\end{equation}
By considering the constant observable $\bm h(\bm x) = 1$ and substituting into Eq.~\eqref{eq:WBA}, it is immediately clear that the sum of the coefficients must satisfy
\begin{equation}
\label{eq:sum-to-one}
    \mathcal{BR}[1](\bm x) = \sum_{j=-J}^J c_j = 1.
\end{equation}
In addition, because dynamics on tori are conservative, we enforce the time-reversible ``palindromic'' constraint
\begin{equation}
\label{eq:palindromic}
    c_{-j} = c_j\quad \text{for all} \quad 1\leq j \leq J.    
\end{equation}
There are $J+1$ independent constraints in Eqs.~\eqref{eq:sum-to-one} and \eqref{eq:palindromic}, leaving $J$ unknowns which fully determine $\bm c$.

Because the ergodic average is unique on a trajectory, it must be the case that for all $t\in\Zbb^+$
\begin{equation*}
    \mathcal{B}[\bm h](\sympmap^{t+1}(\bm x)) = \mathcal{B}[\bm h](\sympmap^t(\bm x)).
\end{equation*}
For $DT\geq J$ (recall $D$ is the output dimension of $\bm h$), Birkhoff RRE uses this fact to define a residual
\begin{equation*}
    R_{\mathrm{RRE}}^2 = \frac{1}{T C_{\mathrm{RRE}}} \sum_{t=0}^{T-1}\norm{\mathcal{BR}[\bm h](\sympmap^{t+1}(\bm x))-\mathcal{BR}[\bm h](\sympmap^t(\bm x))}^2,
\end{equation*}
where the normalization constant $C_{\mathrm{RRE}}$ is
\begin{equation*}
    C_{\mathrm{RRE}} = \sum_{t=0}^{T+2J-1} w_t \norm{\bm h \circ F^{t+1}(x)-\bm h\circ \sympmap^t(x)}^2 = \mathcal{WB}_{T+2J}\left[\norm{\bm h \circ \sympmap - \bm h}^2\right](x).
\end{equation*}
The residual $R_{\mathrm{RRE}}$ can be thought of as measuring the consistency of the extrapolated Birkhoff average. 
By making the average independent of the initial point of the trajectory, it makes it so that $\mathcal{BR}[\bm h](\bm x)$ goes to the true value as $J$ increases.
Substituting Eq.~\eqref{eq:WBA}, the residual can be written to explicitly include $\bm c$ as
\begin{equation*}
    R_{\mathrm{RRE}}^2 = \frac{1}{T C_{\mathrm{RRE}}}\sum_{t=0}^{T-1} \norm{\sum_{j=-J}^J c_j \bm u_{j+J+t}}^2, \qquad C_{\mathrm{RRE}} = \sum_{t=0}^{2J+T-1} w_t \norm{\bm u_t}^2, \qquad \bm u_{t} = \bm h(\sympmap^{t+1}(\bm x)) - \bm h(\sympmap^t(\bm x)).
\end{equation*}
This explicitly shows that the residual is quadratic in the unknowns, demonstrating that it can be equaivalently written in matrix form as
\begin{equation*}
    R_{\mathrm{RRE}}^2 = \frac{1}{T C_{\mathrm{RRE}}} \bm c^T U^T U \bm c, \qquad U = \begin{pmatrix}
        \bm u_0     & \dots  & \bm u_{2J} \\
        \vdots      & \ddots & \vdots \\
        \bm u_{T-1} & \dots  & \bm u_{2J+T-1}
    \end{pmatrix}.
\end{equation*}
So, we find that the requirement that $D T\geq J$ is equivalent to the possibility of $U$ being at least rank $J$.
In total, the least-squares problem for the coefficients $\bm c$ is
\begin{gather}
\label{eq:RRE}
    R_{\mathrm{RRE}}^2 = \min_{\bm c} \frac{1}{T C_{\mathrm{RRE}}} \bm c^T U^T U \bm c, \\ 
    \text{s.t. } \sum_{j=-J}^J c_j = 1, \qquad c_{-j} = c_j \quad\text{for all } 1\leq j \leq J.
\end{gather}
After substituting the constraints into the optimization, this is a least-squares problem for $J$ unknowns that can be solved QR factorization (although, one must be careful in order to achieve good floating point behavior, see Appendix \ref{app:solving-RRE}).
It is simple to observe that $R_{\mathrm{RRE}}$ for an invariant torus must go to zero at least as fast as a weighted Birkhoff average.
Moreover, for $d=1$, the rate of convergence of the residual was proven to be faster than the best known rate for weighted Birkhoff averages \cite{ruth_finding_2024}.
The convergence theory for higher dimensions has not been investigated, however.

% In practice, we solve \eqref{eq:RRE} adaptively.
% This is performed by choosing a desired residual $\epsilon$, fixing $T$, and increasing $J$ iteratively until $R_{\mathrm{RRE}} < \epsilon$.
% These updates can be performed by QR column updates as in App.~\ref{app:qr-updates}, so the worst-case complexity is still $\Ocal(TD J^2)$.

Assuming $c_J \neq 0$, Eq.~\eqref{eq:WBA} can interpreted as a linear difference equation to predict the evolution of the observed trajectory as
\begin{equation}
\label{eq:extrapolation}
    \bm h(\sympmap^{2J}(\bm x)) - \bm h_0 = - \frac{1}{c_J} \sum_{j = 0}^{2J-1} c_j (\bm h(\sympmap^{j}(\bm x))-\bm h_0).
\end{equation}
This is the origin of the word ``extrapolation'' in Birkhoff RRE.
From the point of view of linear difference equations, it is straightforward to see how to obtain the trajectory's frequencies once $\bm c$ is known.
Like linear homogeneous differential equations, linear homogeneous difference equations of the form \eqref{eq:extrapolation} have $2J$ homogeneous solutions $\varphi_j : \Zbb \to \Cbb$ with eigenvalues $\psi_j\in\Cbb$ for $1\leq \abs{j} \leq J$ satisfying
\begin{equation*}
    \varphi_j(n) = \psi_j^n, \qquad \sum_{t=-J}^J c_t \varphi_j(t) = 0.
\end{equation*}
Because each $c_j$ is real and we have enforced the palindromic constraint \eqref{eq:palindromic}, it must be the case that if $\psi_j\in\Cbb$ is an eigenvalue, its complex conjugate $\psi_{-j} = \overline \psi_j$, inverse $\psi_j^{-1}$, and conjugate inverse $\overline \psi_j^{-1}$ are all also eigenvalues (this is the standard eigenvalue symmetry of a Hamiltonian matrix).
As a result, most eigenvalues $\psi_j$ of the optimal filter are fixed to the unit circle $\abs{\psi_j}=1$, where the inverse and conjugate coincide.
Assuming there are $J^*$ of these eigenvalues where $1 \leq J^* \leq J$, we define the frequencies $\Omega_j \in \Tbb$ as
\begin{equation*}
    \Omega_j = \arg \psi_j := \frac{1}{2\pi \I} \log \psi_j \mod 1.
\end{equation*}
The eigenvalues are simultaneously computed by a single $J\times J$ eigenvalue problem in $\Ocal(J^3)$ time; see \cite{ruth_finding_2024} for details.

Finally, once the frequencies are known, they can be projected back onto the signal by a single least-squares problem
\begin{equation*}
    \min_{\{\bm H_j\}} \sum_{t = 0}^{T+2J-1} \norm{\bm h(\sympmap^t(\bm x)) - \sum_{j=-J^*}^{J^*} \bm H_j e^{2\pi \I \Omega_j t}},
\end{equation*}
where we use the convention that $\Omega_0 = 0$.
This, like the coefficients $\bm c$, is computed by QR factorization in $\Ocal(TD J^2)$ time, leading to an overall complexity of $\Ocal(TD J^2)$.
We denote the squared magnitude of the coefficients as $H_j^2 = \norm{\bm H_j}^2$, and sort the pairs $(\Omega_j, H_j^2)$ in order of the magnitudes so that $H_1^2 \geq H_2^2 \geq \dots \geq H_{J^*}^2$.

There is currently no theory for the convergence of the frequencies $\Omega_j$ to true values $\bm \omega \cdot \bm k_j$ as $J\to \infty$.
However, as has already been indirectly demonstrated by the quality of validation error in Sec.~\ref{sec:parameterization}, we find that the frequencies $\Omega_j$ often converge to near floating point accuracies in practice.
This is particularly true for more prominent modes, corresponding to smaller values of $j$. 
As $j$ increases we find the approximations to $\Omega_j$ degrade.
We conjecture the rate of convergence of individual frequencies of analytic functions is $\Ocal(e^{-r J^{1/d}})$ for some $r > 0$, matching the rate of convergence of Fourier coefficients.

\subsection{Finding an initial rotation vector}
\label{subsec:initial-rotation}
Now, we would like to determine the underlying rotation vector and wavenumbers.
The input for this subsection is the output of the previous, namely a frequency vector $\bm \Omega = (\Omega_1, \dots, \Omega_{\JMAP})$ with corresponding the magnitudes $\bm H = (H_1^2, \dots, H_{\JMAP}^2)$, where $\JMAP$ is the number of frequencies that we will use for this step, satisfying $1\leq\JMAP\leq J^*$.
The value of $\JMAP$ is defined by the user.
As $\JMAP$ increases, so does the cost of finding the rotation vector, so it is practically important for it to be smaller than $J^*$.
For numerical stability, it is also important that $\JMAP$ not include frequencies where $H_{\JMAP} \ll H_1$, so always choose $J_0$ so that $H_{\JMAP}/H_1 > \epsilon_{\mathrm{MAP}}$ (where the default value is $\epsilon_{\mathrm{MAP}} = 10^{-3}$).
Using the frequencies and magnitudes, we estimate the rotation vector and wavenumbers according to the MAP estimation problem \eqref{eq:map}, where we maximize the posterior $\prob(\bm \omega, \bm \pi \mid \bm \Omega, \bm H)$ over candidate rotation vectors $\bm \omega \in \Tbb^d$ and wavenumber mappings $\bm \pi : [L]\to\Zbb^d$ from the index $j$ to wavenumber $\bm k$.
The purpose of MAP estimation is to encode our belief for the most likely rotation vector and wavenumbers conditioned on our data.
The priors for this estimation will include both a smoothness assumption on the torus coefficients $\bm H$ and a basic model for the accuracy of the frequencies $\bm \Omega$ returned by Birkhoff RRE.

To see why the prior of a smooth torus is necessary for the problem, consider a candidate rotation vector $\bm \omega$ with at least one irrational element $\omega \in \bm \omega$.
Because $\omega$ is irrational, integer multiples of the frequency are dense in $\Tbb$; i.e.~for any element of the frequency data $\Omega_j$, there exists a sequence of wavenumbers $k_\ell\in\Zbb$ such that $\Tnorm{\omega k_\ell -\Omega_j} \to 0$ as $\ell \to \infty$.
This is problematic because any errors in the frequency measurements can be pathologically explained by an unbounded sequence of wavenumbers.
The smooth-torus prior is meant to counteract this type of overfitting by appropriately weighting the likelihood of observing a wavenumber by the probability that the associated mode is prominent in the torus.
Additionally, the probabilistic framework helps decide between wavenumbers when multiple wavenumbers could explain a frequency due to the problem of small denominators.
In this way, while the underlying problem is not fundamentally probabilistic, the MAP formalism gives a straightforward method to analyze the frequency output of the previous step.

We start by defining a prior distribution on the invariant torus of
\begin{equation*}
    \bm h(S(\bm \theta)) = \sum_{\bm k} \bm h_{\bm k} e^{2\pi \I \bm k \cdot \bm \theta}, \qquad \bm h_{\bm k} =\overline{\bm h}_{-k}, \qquad \real(\bm h_{\bm k}), \imag(\bm h_{\bm k}) \sim \mathcal N\left(0, \frac{\sigma(\bm k)^2}{2} I\right),
\end{equation*}
where $\sigma(\bm k)$ for $\bm k\in\Zbb^d$ are wavenumber-dependent hyperparameters and $\mathcal{N}(\bm \mu, \Sigma)$ is a joint-normal distribution with mean $\bm \mu$ and covariance $\Sigma$ defined by the probability distribution function (PDF) for $\bm x \sim \mathcal{N}(\bm \mu, \Sigma)$ as
\begin{equation*}
    \prob(\bm x) = \frac{1}{\sqrt{2 \pi \det \Sigma}}\exp\left(-\frac{1}{2} (\bm x - \bm \mu)^T\Sigma^{-1} (\bm x - \bm \mu) \right).
\end{equation*}
The variances are specifically chosen so that $\bm h$ is analytic by taking
\begin{equation*}
    \sigma(\bm k) = C e^{-r \norm{\bm k}},
\end{equation*}
where $C,r>0$ are hyperparameters (another choice could be $\sigma(\bm k) = C \norm{\bm k}^{-N}$ for $C^N$ tori, but the analytic case appears to be typical for smooth dynamical systems).
To find these parameters, consider the idealized case that $\bm H$ consist of a sorted list of the exact expected norms $H_k^2 = \expec[\norm{\bm h_{\bm k}}^2] = D\sigma(\bm k)^2$. 
In such a case, the list is sorted according to increasing wavenumbers in the lattice, which is asymptotically given by $\norm{\bm \pi(j)} \sim ((j-1)/V_d)^{1/d}$ where $V_d$ is the unit volume of a sphere in $d$ dimensions. 
Using this assumption, we heuristically fit a line to the log-norm of the coefficients $\bm H$
\begin{equation}
\label{eq:hp-optimization}
    \min_{C, r} \sum_{j=1}^{J_{C,r}} \abs{\log{\norm{\bm H}_j^2} - \log\left(D \sigma\left(\left((j-1)/V_d\right)^{1/d}\right)^2\right)}^2,
\end{equation}
When this problem is expanded, it becomes a simple linear regression for $\log C$ and $r$, so the hyperparameters can be quickly obtained from data.
We choose the number of points $J_{C,r}$ to be approximately a third of the total number of available values from the output of the extrapolation step.

Because the observed coefficients are observed in sorted norm, we will use the distribution of the squared magnitude of the coefficients $\norm{\bm h_{\bm \pi(j)}}^2 \sim \chi^2(\sigma(\bm \pi(j))^2, D)$ with the probability density function
\begin{equation}
\label{eq:Prob_Hk}
    \prob(\norm{\bm h_{\bm k}}^2=x^2) = \frac{1}{(2\sigma(\bm k)^2)^{D/2}\Gamma(D/2)}x^{D - 2} \exp\left(-\frac{x^2}{2\sigma(\bm k)^2}\right).
\end{equation}
We assume that this probability is independent of the rotation vector, i.e.~$\prob(\norm{\bm h_{\bm k}}^2 \mid \bm \omega) = \prob(\norm{\bm h_{\bm k}}^2)$.

The extrapolation step does not return wavenumbers, so we cannot directly use the probability \eqref{eq:Prob_Hk} in our optimization.
So, to link these probabilities with the observed data $\bm H$, we assume that extrapolation step returns an order statistic $\bm \pi : [\JMAP]\to \Zbb^d$ on the magnitudes $H_j^2 = \norm{\bm h_{\bm \pi(j)}}^2$ with no measurement error. 
An order statistic is defined by the requirement that the magnitudes are sorted as $H_1^2 \geq H_2^2 \geq \dots \geq H_{\JMAP}^2$ and $\norm{\bm h_{\bm k}}^2 \leq H_{\JMAP}^2$ for all unobserved $\bm k \notin \bm \pi([\JMAP])$.
We note that order statistics on non-identically independently distributed random variables are expensive to compute, requiring the factorial cost of computing permanents \cite{david_order_2003}.
However, we will find that we only require the joint probabilities $\prob(\bm H, \bm \pi)$, allowing us to avoid this particular computational complexity.

In particular, the joint probability of observing $(\bm H,\bm \pi)$ is
\begin{align}
\label{eq:Hpi-joint}
    \prob(\bm H, \bm \pi) &= \Delta_{\bm\pi} \left(\prod_{j=1}^{\JMAP}\prob(H_{\bm \pi(j)}^2=H_j^2) \right)\left( \prod_{\bm k\in\Zbb^d\backslash\bm (\bm \pi([\JMAP])\cup\{0\})}\prob(H_{\bm k}^2 \leq H_{\JMAP}^2)\right), \\
\label{eq:Hpi-joint2}
    &\propto \Delta_{\bm \pi} \prod_{j=1}^{\JMAP} \frac{\prob(H_{\bm \pi(j)}^2=H_j^2)}{\prob(H_{\bm \pi(j)}^2 \leq H_{\JMAP}^2)},
\end{align}
where 
\begin{align*}
    \Delta_{\bm \pi} %&= \prod_{j=1}^{\JMAP}\prod_{j'=j+1}^{\JMAP}\Delta_{\bm \pi(j),\bm\pi(j')}, \\
    &= \begin{cases}
        1, & \bm \pi(j)\neq\bm\pi(j') \quad \text{ for all } \quad 1\leq j<j'\leq \JMAP,\\
        0, & \text{else}.
    \end{cases}
\end{align*}
In words, the formula for the joint probability \eqref{eq:Hpi-joint} is the probability that each individual entry of $H_j \in \bm H$ is observed according to the distribution defined by $\bm \pi(j)$, multiplied against the probability that every other Fourier coefficient has a magnitude less than the final entry of $H_{\JMAP}$.
The function $\Delta_{\bm \pi}$ disallows repeated draws of the same wavenumber.
We note that the proportionality constant in \eqref{eq:Hpi-joint2} only depends on the data, and not on $\bm \pi$.

For the rotation vector, we assume a uniform prior on $\Tbb^d$ independent of the coefficients $\bm \omega \sim \mathcal{U}(\Tbb^d)$, corresponding to a probability distribution function
\begin{equation}
\label{eq:omega-prior}
    \prob(\bm \omega) = 1.
\end{equation}
The measurements of frequencies $\Omega_j$ are assumed to be sampled from a circular normal with mean $\bm \omega \cdot \bm \pi(j)$ and standard deviation $\sigma_\omega$
\begin{equation*}
    \prob(\Omega_j \mid \bm \omega, \bm \pi) = \frac{1}{C} \exp\left(\frac{\sin^2((\Omega_j - \bm \omega \cdot \bm \pi(j))\pi)}{2 \pi^2 \sigma_\omega^2}\right),
\end{equation*}
where $C$ is a normalizing constant.
In the limit as $\sigma_\omega \to 0$, this distribution approaches
\begin{equation*}
    \prob(\Omega_j \mid \bm \omega, \bm \pi) = \frac{1}{\sqrt{2\pi \sigma_\omega^2}} \exp\left(\frac{\Tnorm{\Omega_j - \bm \omega \cdot \bm \pi(j)}^2}{2 \sigma_\omega^2}\right),
\end{equation*}
which is what we use in practice. 
The hyperparameter $\sigma_\omega$ is typically chosen to be very small (i.e.~$\Ocal(10^{-7})$ to $\Ocal(10^{-11})$), and can be adjusted to the convergence of the extrapolation step. 
In particular, when the map is known to be numerically symplectic --- either from an analytic form or a symplectic integrator --- the tolerance can be chosen much smaller than if, e.g., the map is obtained via a Runge-Kutta method.
We assume that the frequencies $\Omega_j$ are independent of $\bm H$ and each other, so that
\begin{equation*}
    \prob(\bm \Omega \mid \bm \omega, \bm \pi) = \prob(\bm \Omega\mid \bm H, \bm \omega, \bm \pi) = \prod_{j=1}^{\JMAP}\prob(\Omega_j \mid \bm \omega, \bm \pi(j))
\end{equation*}

Now that we have defined our priors, we return to the Bayesian MAP estimation problem
\begin{equation*}
    \max_{\bm \omega, \bm \pi} \prob(\bm \omega, \bm \pi \mid \bm \Omega, \bm H),
\end{equation*}
where we use the shorthand $\bm \Omega = (\Omega_1, \dots, \Omega_{\JMAP})$ and $\bm H = (H_1, \dots, H_{{\JMAP}})$ for the first ${\JMAP}$ frequency and coefficient measurements.
Applying Bayes' rule, we have
\begin{align*}
    \prob(\bm \omega,\bm \pi \mid \bm \Omega, \bm H) &= \frac{\prob(\bm \omega)}{\prob(\bm \Omega, \bm H)} \prob(\bm \Omega, \bm H, \bm \pi \mid \bm \omega), \\
    &\propto \prob(\bm \Omega, \bm H, \bm \pi \mid \bm \omega),
\end{align*}
where we used Eq.~\ref{eq:omega-prior} and the fact that the marginal probability $\prob(\bm \Omega, \bm H)$ is independent of the optimization variables to put it into the proportionality constant.
Because we have assumed that $\bm \Omega$ is independent of $\bm H$ and that $\bm H$ is independent of $\bm \omega$, this reduces to
\begin{align}
\nonumber
    \prob(\bm \omega, \bm \pi \mid \bm \Omega, \bm H) &\propto \prob(\bm \Omega \mid \bm \omega, \bm \pi)\prob(\bm H, \bm \pi), \\
\label{eq:map-explicit}
    &\propto \Delta_{\bm \pi} \prod_{j=1}^{\JMAP} \alpha_j(H_j,\bm\pi(j))\prob(\Omega_j \mid \bm \omega, \bm \pi(j)), && \alpha_j(H_j,\bm \pi(j)) = \frac{\prob(H_{\bm \pi(j)}^2=H_j^2)}{\prob(H_{\bm \pi(j)}^2 \leq H_{\JMAP}^2)}.
\end{align} 

Eq.~\eqref{eq:map-explicit} is the desired outcome of a weighted probability of observing the frequencies $\bm \Omega$.
The frequency probabilities penalize deviations between the candidate rotation vector and the observed frequencies, which is our primary objective.
These probabilities are then weighted by the likelihood that the torus coefficient magnitude deviates too much from the analytic smooth-torus assumption.
In the case that $\Delta_{\bm \pi} = 1$, the logarithm of the proportional value of the posterior 
\begin{equation}
\label{eq:Lmap}
    L_{\mathrm{MAP}} = \sum_{j=1}^{\JMAP} \log \alpha_j + \log \prob(\Omega_j \mid \bm \omega, \bm \pi(j))
\end{equation}
can be used as an approximate measure of fit.
However, this measure can be unreliable, as we have found that correct rotation vectors can sometimes have large negative values of $L_{\mathrm{MAP}}$.

% If the coefficients satisfied $\alpha_j = 1$, the sum in Eq.~\eqref{eq:map} would have the pathological problem of arbitrarily large wavenumber solutions. 
% Instead, they exponentially penalize large values of $\bm \pi(j)$, as $\prob(H_{\bm \pi(j)}^2=H_j^2) \to 0$ and $\prob(H_{\bm \pi(j)}^2 \leq H_{\JMAP}^2) \to 1$ as $\norm{\bm \pi(j)}\to\infty$.

Now, we consider the process of solving the MAP problem \eqref{eq:map}.
Before we describe the method, first note that the problem is discrete in $\bm \pi$, indicating an optimization that is potentially very expensive. 
So, we will not perform a full global optimization

First, instead of optimizing over all possible values of $\bm \omega$, we will restrict to values such that 
\begin{equation*}
    \bm \omega = \begin{pmatrix}
        \Omega_{j_1} \\ \vdots \\ \Omega_{j_d}
    \end{pmatrix}\subset \bm \Omega.
\end{equation*}
This assumes that there exists a valid rotation vector comprised of the frequencies that we have observed from the RRE process.
From this, we have reduced the optimization problem to
\begin{equation*}
    \min_{\bm \omega \subset \bm \Omega} \min_{\bm \pi} \prob(\bm \omega, \bm \pi \mid \bm \Omega, \bm H).
\end{equation*}

Then, for the wavenumber mapping $\bm \pi$, we first restrict to finite number of wavenumbers in the finite infinity-norm grid
\begin{equation*}
    \norm{\bm \pi(j)}_\infty = \max_{\ell \in [d]}\abs{(\pi(j))_\ell} \leq P, \qquad 1\leq j \leq \JMAP, 
\end{equation*}
for some fixed value of $P>0$.  
By default, we consider $P$ to scale as $\gamma (\JMAP)^{1/d}$, where $\gamma \approx 10$. 
Then, on that grid, we greedily choose the values $\bm \pi(j)$ by minimizing
\begin{equation}
\label{eq:MAP-optimization}
    \min_{\substack{\norm{\bm \pi(j)}_\infty \leq P \\ \bm \pi(j) \notin \{\bm \pi(1), \dots, \bm \pi(j-1)\}}} \alpha_j(H_j, \bm \pi(j)) \prob(\Omega_j \mid \bm \omega, \bm \pi(j)).
\end{equation}
Note that by excluding the previous wavenumber draws in the minimization, we automatically ensure that $\Delta_{\bm \pi} = 1$. 
If the result would not change if we allowed for repeated draws, the greedy optimization gives the true optimal posterior conditioned with fixed rotation vector on the allowed grid.

The primary computational cost of the MAP step is in the optimization \eqref{eq:MAP-optimization}.
A direct implementation requires searching over $(2P+1)^d$ potential matches for all $1\leq j \leq \JMAP$, resulting in a total cost that scales as $\Ocal(\JMAP P^d)$.
To accelerate the optimization, we first sort the allowed wavenumbers $A = \{\bm k | \norm{\bm k}_\infty \leq P \}$ by the frequencies $\bm \omega \cdot \bm k$.
Because the probability $\prob(\Omega_j \mid \bm \omega, \bm \pi(j))$ is sharply localized about the frequency $\Omega_j$, we only have to search for $\bm \pi(j)$ in the sorted list over a small fixed number of frequencies about $\Omega$ in $A$ (approximately 5 in each direction).
This reduces the cost of computing wavenumbers for fixed $\bm \omega$ to $\Ocal(P^d \log P^d)$ --- the time to sort the list.
Iterating over all allowed rotation vectors within $\bm \Omega$, the total cost is then $\Ocal((\JMAP)^d P^d \log P^d)$, which scales as $\Ocal((\JMAP)^{d+1} \log \JMAP)$ when $P \propto (\JMAP)^{1/d}$.
We note that there are cases where $P$ scaling as $\JMAP$ is inappropriate --- e.g.~for very thin tori where most of the motion is occuring in a single dimension. 
In such cases one might consider taking $P$ larger, or perhaps choosing an anisotropic grid to search over wavenumbers.

\subsection{Refining the rotation vector}
\label{subsec:refine-rotation}
If the MAP estimation succeeded, the output from Sec.~\ref{subsec:initial-rotation} will be a valid rotation vector.
However, while this rotation vector is valid, we have not necessarily addressed the issue of whether the found rotation vector is the best one for Fourier approximation. 

To determine a precise notion of ``best rotation vector,'' we consider a geometric approach. 
As discussed in \ref{subsec:nonuniqueness}, when we transform our rotation vector, we equivalently transform the homology generators of the torus, meaning we need some metric of what the best homology generators are.
To do so, consider $S : \Tbb^d \to \Rbb^D$ to be the invariant torus in Euclidean space for a given rotation vector.
A loop around the torus $S$ in the $\theta_1$ direction is parameterized by fixing $\theta_m$ for $1 < m \leq d$ and varying $0\leq \theta_1 < 1$. 
Then, the length of that loop is given by
\begin{equation*}
    L_1(\theta_2, \dots, \theta_d) = \int_{\Tbb}\left(\pd{(\bm h \circ S)}{\theta_1}\right)^T \left(\pd{(\bm h \circ S)}{\theta_1}\right) \dif \theta_1.
\end{equation*}
Our goal is to find the direction $\theta_1$ with a shortest average loop that maintains the action-angle conjugacy in the $L^2$ sense, which is defined for a fixed torus as
\begin{align*}
    \overline L_1 &= \int_{\Tbb^d}\left(\pd{(\bm h \circ S)}{\theta_1}\right)^T \left(\pd{(\bm h \circ S)}{\theta_1}\right) \dif \bm \theta, \\
    &= \bm e_1^T \left[ \int_{\Tbb^d}\left(\pd{(\bm h \circ S)}{\bm \theta}\right)^T \left(\pd{(\bm h \circ S)}{\bm \theta}\right) \dif \bm \theta \right]\bm e_1, \\
    &= \bm e_1^T G \bm e_1.
\end{align*}
The matrix $G$ can be interpreted as the averaged metric of the observable embedding of the invariant torus in $\Rbb^d$. 

When we consider a transformation of the invariant torus with $\bm \theta' = A \bm \theta$ for $A \in GL(d,\Zbb)$, the metric changes as
\begin{align*}
    G' &= \int_{\Tbb^d}\left(\pd{\bm \theta}{\bm \theta'}\right)^T\left(\pd{(\bm h \circ S)}{\bm \theta}\right)^T \pd{(\bm h \circ S)}{\bm \theta}\pd{\bm \theta}{\bm \theta'} \dif \bm \theta, \\
    &= A^{-T} G A^{-1}.
\end{align*}
A shortest average loop direction around the torus is then defined by the optimization
\begin{equation*}
    \overline L_1^\star = \min_{A \in GL(d,\Zbb)}\bm e_1^T (A^{-T} G A^{-1}) \bm e_1.
\end{equation*}
(Note that this minimization does not have a unique minimizer even in the 1D case, as $A \to -A$ produces the same solution.)
Let $B$ be any matrix square root of the average metric $G = B^T B$ (e.g.~the Cholesky decomposition). 
We can rewrite the optimization problem as
\begin{equation*}
    \overline L_1^\star = \min_{A \in GL(d,\Zbb)}\norm{B A^{-1} \bm e_1}^2.
\end{equation*}
When written this way, it becomes apparent that a shortest average loop is equivalent to the Shortest Vector Problem (SVP), which is a well known problem in lattice cryptography 
\cite{micciancio_complexity_2002}. 
This problem is NP-hard for randomized reductions in the dimension $d$ \cite{ajtai_generating_1996}, meaning the fastest algorithms for computing $L_1^\star$ are believed to be very expensive as $d$ grows.
In the case of large dimension, a well known approximation algorithm for SVP is the Lenstra-Lenstra-Lov\'asz (LLL) algorithm \cite{lenstra_factoring_1982}. 
However, in the case of the invariant tori we consider, $d$ is typically very small, meaning the exact problem can be solved quickly.

The shortest vector can always be used as the first column of the transformed matrix $B'$. 
From here, however, there are multiple ways this can generalize to higher dimensions.
For this purpose, we choose the Korkine-Zolatarev (KZ) basis \cite{micciancio_complexity_2002}, which is defined so that the transformed $B$ has short and nearly orthogonal columns.

To define the KZ basis, let $\bm b_n$ be the $n$th column of $B$.
The orthogonal lattice Gram-Schmidt basis is defined by 
\begin{equation*}
    \hat{\bm b}_n = \bm b_n - \sum_{m=1}^{n-1} \mu_{mn} \hat{\bm b}_m, \qquad \mu_{mn} = \frac{\bm b_n \cdot \hat{\bm b}_m }{|\hat{\bm b}_m|^2}.
\end{equation*}
Then, we define the $n$th projected lattice $\Lambda_n$ as the set of points
\begin{equation*}
    \Lambda_n = \{\bm x \mid \bm x = P_n(B \bm k), \bm k\in\Zbb^d \}, \qquad P_n(\bm x) = \bm x - \sum_{m=1}^{n-1} \frac{\bm x \cdot \hat{\bm b}_m }{|\hat{\bm b}_m|^2}\hat{\bm b}_m.
\end{equation*}
Note that $\Lambda_1$ is just the original lattice defined by $B$.
We define a basis $B$ as KZ reduced if 
\begin{itemize}
    \item $\hat{\bm b}_n$ is a shortest vector in $\Lambda_n$, and
    \item for all $n>m$, $\abs{\mu_{mn}} \leq 1/2$.
\end{itemize}

In the definition of a KZ basis, the two requirements correspond directly to notions of shortness and orthogonality.
The first requirement exactly coincides with $\hat{\bm b}_1$ being a minimizer of the SVP for the shortest average loop.
The following columns are minimizers of SVPs on toroidal coordinates that are orthogonal to the previously computed directions.
The second requirement can be recognized as the cosine of the angle between two vectors, implying that the angles between any two basis vectors is greater than $\pi/3.$

To appreciate the usefulness of the orthogonality requirement to the parameterization, we note that in Fourier space on the torus, the averaged metric becomes
\begin{equation}
\label{eq:Fourier-G}
    G = \sum_{\bm k \in \Zbb^d} \norm{\bm h_{\bm k}}^2 \bm k \bm k^T .
\end{equation}
In this way, $G$ can be thought of as the second moment of the distribution of Fourier coefficients (because $\norm{\bm h_{\bm k}} = \norm{\bm h_{-\bm k}}$, the first moment must always be zero).
So, the requirement of small angles in the basis is equivalent to the requirement that a Gaussian approximation of the Fourier spectrum is well aligned with the coordinate axes: i.e.~we select the left hand side of Fig.~\ref{fig:Fourier-shear}.

The equation \eqref{eq:Fourier-G} also gives a method to compute the matrix $G$.
To estimate this, we use the Fourier labeling $\bm \pi$ determined in the previous section (possibly for different values of $\JMAP$ and $P$, as this must only be computed once).
Using that mapping, we can approximate $G$ as
\begin{equation*}
    G \approx \sum_{j=1}^{\JMAP}H_j^2 \bm \pi(j) \bm \pi(j)^T.
\end{equation*}
From here, we use an implementation of a KZ reduction algorithm \cite{zhang_hkz_2012} by the LLLplus.jl package.
Once we find a KZ reduced basis, we post-process to restrict the non-uniqueness of the solution.
In particular, if $\bm \omega' = A \bm \omega$ is the transformed rotation vector, our final step is to multiply transform $A$ to $A \Diag(\mathrm{sgn}(\bm \omega'))$, where the sign function $\mathrm{sgn} : \Tbb \to \Zbb$ is defined by
\begin{equation*}
    \mathrm{sgn}(\omega) = \begin{cases}
        1, & 0 \leq \omega \leq 1/2,\\
        -1, & \text{else}.
    \end{cases}
\end{equation*}
The result is that the final rotation vector $\bm \omega'$ is chosen to have entries in $[0.0,0.5] \mod 1$.

\subsection{Algorithm adjustments for islands}
\label{subsec:islands}
While invariant tori are the typical integrable structure of symplectic maps, it is possible that trajectories lie on a periodic sequence of tori that we refer to as an \textit{island chain}.
We define the period of an island chain including $x$ as the minimal value $p \in \Zbb^+$ such that $x$ is on an invariant torus of $\sympmap^p$ (in this way, an invariant torus is an island chain with period $1$).
Then, we define an island chain with period $p$ as a sequence of tori $S^{(j)} : \Tbb^d \to \mathcal{M}$ for $1\leq j \leq p$ with rotation vector $\bm \omega$ that satisfies
\begin{equation*}
    \sympmap \circ S^{(j)} = \begin{cases}
        S^{(j+1)}, & 1 \leq j < p, \\
        S^{(1)} \circ \tau_{\bm \omega}, & j = p.
    \end{cases}
\end{equation*}

When transformed to the frequency domain, the signal $\bm h \circ \sympmap^t$ has the form
\begin{equation*}
    \bm h \circ \sympmap^t(x) = \sum_{j=0}^{p-1} \sum_{\bm k \in \Zbb^d} \bm h_{j\bm k} \lambda_{j\bm k}^t, \qquad \lambda_{j\bm k} = e^{2\pi \I (j + \bm \omega \cdot \bm k)/p}.
\end{equation*}
So, from the perspective of the signal, an island chain can be equivalently be thought of as a dimension $d+1$ invariant torus with a resonant rotation vector $\bm \omega_{\mathrm{isl}} = (1/p,\bm \omega/p)$.
We can use this insight to easily generalize the algorithm above. 

To find $p$, we use the same method as \cite{ruth_finding_2024}. 
Let $\epsilon_{\mathrm{isl}} > 0$ be a tolerance and $p_{\mathrm{max}}$ be the maximum value of $p$.
Then, we identify a trajectory as having a period $p$ if there is any element of the frequency vector $\Omega_j \in \bm \Omega$ which satisfies
\begin{equation*}
    \Tnorm{\Omega_j - \frac{n}{p}} < \epsilon_{\mathrm{isl}}, \qquad n,p\in \Zbb, \qquad 1\leq p\leq p_{\mathrm{max}}.
\end{equation*}
Typical values of these constants are $p_{\mathrm{max}} = 10$ and $\epsilon_{\mathrm{ada}} = 10^{-8}$.

If the torus is recognized as an island, the MAP estimation is essentially the same, except we would like to learn $\bm \pi : [\JMAP] \to [p] \times \Zbb^d$, where the first element corresponds to the discrete frequency associated with the rational frequency $1/p$.
For the prior on $\bm h_{j \bm k}$, we use
\begin{equation*}
    \real \bm h_{j\bm k}, \imag \bm h_{j\bm k} \sim \mathcal{N}\left(0, \frac{\sigma(\bm k)^2}{2} I\right),
\end{equation*}
i.e.~the frequency $j$ is not included the prior. 
Because the period $p$ is determined before the MAP estimation, it is treated as a known quantity. 
This means that the prior for $\bm \omega$ is not changed.

Finally, for rotation vector refinement, we perform KZ reduction on the averaged metric
\begin{equation*}
    G_{mn} = \sum_{j=1}^p \int_{\Tbb^d} \pd{(\bm h \circ S^{(j)})}{\theta_m} \cdot \pd{(\bm h \circ S^{(j)})}{\theta_n} \dif \bm \theta.
\end{equation*}
The diagonal components $G_{mm}$ can be interpreted as the sum-of-squares of the average loop length of each torus in the island chain.
When transformed to Fourier space, we represent this matrix as
\begin{equation*}
    G = \sum_{\bm k \in \Zbb^d} \left(\sum_{j=1}^p \norm{\bm h_{j\bm k}}^2 \right) \bm k \bm k^T.
\end{equation*}
This can be approximated using $\bm \pi$ from the previous step.

Finally, for the parameterization algorithm in Sec.~\ref{subsec:fourier-adaptive}, we can find parameterizations for each island in the island chain by working directly with the extended frequency $\bm \omega_{\mathrm{isl}}$.
The only adjustment beyond this is that refinement in the rational direction never needs to be adjusted.
That is, we optimize over a resolution $\bm K_{\mathrm{isl}} = (p,\bm K)$, where only $\bm K$ is incremented.

\section{Examples}
\label{sec:examples}
Now that the method has been introduced, we would like to see how it performs under more extensive testing.
We do this in two settings.
In the first example (the coupled standard map, Sec.~\ref{subsec:sm}), we analyze the method's robustness.
We do so by attempting to compute invariant tori for 1000 randomly initialized trajectories and measuring the accuracy of the solutions.
To make the test rigorous, accuracy is computed by an \textit{a posteriori} test of a KAM-like residual, rather than any of the residuals used to compute the torus.
In the second example (the cislunar ER3BP, Sec.~\ref{subsec:cislunar}), we show that the method generalizes to 3D tori in the out-of-plane ER3BP.
The implementation of the algorithm used to perform these examples can be found in the \texttt{SymplecticMapTools.jl} Julia package \cite{Ruth:2025:SymplecticMapTools}.
All reported timings were performed on an Apple M2 CPU.

\subsection{The coupled standard map}
\label{subsec:sm}
We return to the standard map \eqref{eq:coupled-sm} with
\begin{equation*}
    K_{\mathrm{sm}} = \begin{pmatrix}
            0.1 & 0.05 \\ 0.05 & 0.05
        \end{pmatrix},
\end{equation*}
where we use the same observable \eqref{eq:sm-observable}.
The forcing matrix $K_{\mathrm{sm}}$ was chosen to be weak so that there are many invariant tori within the coupled standard map equation. 
The strong off-diagonal component was chosen to be relatively strong, so that there is significant coupling between the variables.
The initial conditions are chosen uniformly at random on the hypercube $(\bm x_0, \bm y_0) \in \Tbb^2 \times [0,1]^2$.

\begin{figure}
    \centering
    \includegraphics[width=1.0\linewidth]{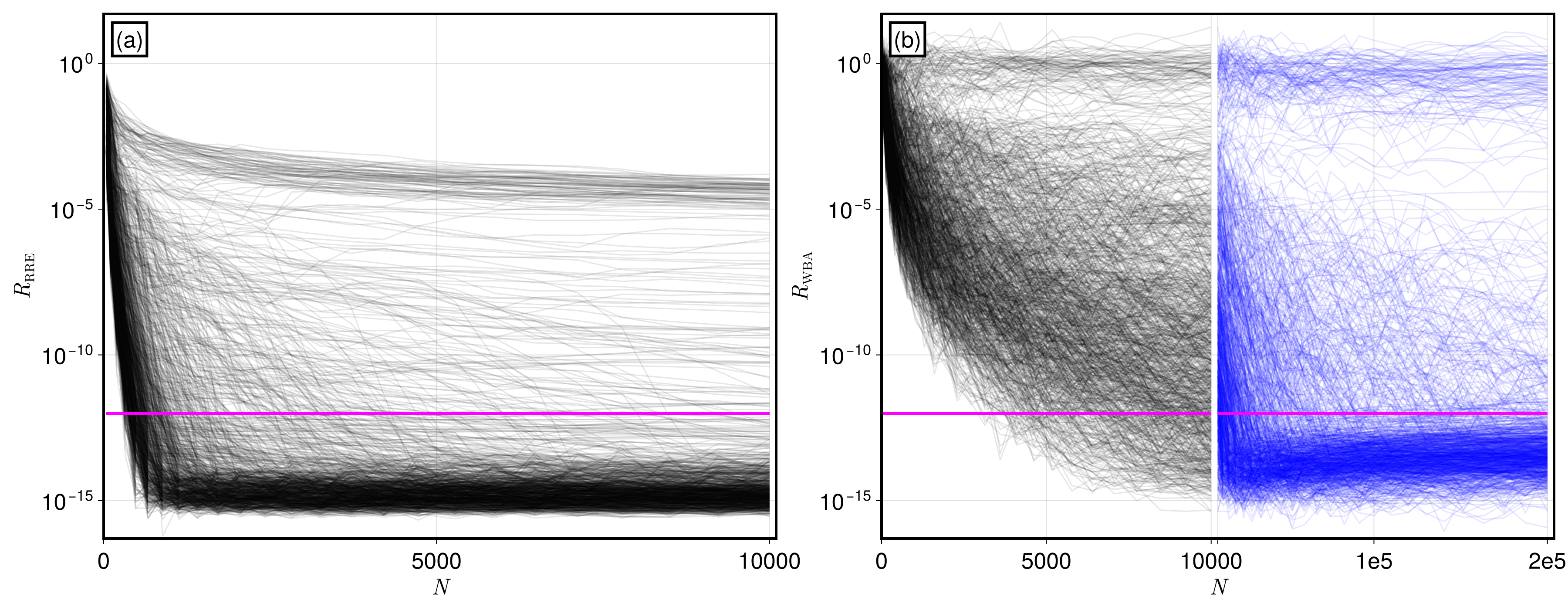}
    \caption{(Left) The Birkhoff RRE residual \eqref{eq:RRE} as a function of trajectory length $N$ for the coupled standard map. 
    The classification tolerance of $1 \times 10^{-12}$ is plotted in magenta. (Right) The Weighted Birkhoff Average residual as a function of trajectory length. 
    The $N$ scaling is broken so that the left half agrees with the Birkhoff RRE horizontal axis and the right half shows the residual for values of $N$ up to $2 \times 10^5$.}
    \label{fig:trajectory-classification}
\end{figure}

To compute the tori, we use the parameters below (where we recall their basic uses)
\begin{itemize}
    \item $(J,T)\in \{(1000,2000),(1500,3000),(2000,4000)\}$ (except for Fig.~\ref{fig:trajectory-classification}): The dimensions of the Birkhoff RRE least-squares problem. 
    We choose the smallest value of these parameters where $R_{\mathrm{RRE}} < 5\times 10^{-14}$ (or, if it is satisfied nowhere, $(J,T) = (2000,4000)$). 
    The corresponding trajectory length is $N=T+2J+1 \in \{4001,6001,8001\}$.
    \item $\JMAP=30$: The number of frequencies used for finding the rotation vector.
    \item $\sigma_\omega = 10^{-10}$: Prior standard deviation of the frequency measurements $\bm \Omega$. 
    \item $K_{\mathrm{max}}=2000$: Maximum number of Fourier modes used to parameterize the torus.
    \item $p_{\mathrm{max}}=10$, $\epsilon_{\mathrm{isl}} = 10^{-8}$: Island identification parameters.
\end{itemize}
As a general rule, as $J,$ $T,$ $J^*,$ and $K_{\mathrm{max}}$ are the most important parameters to adjust, allowing for more accurate computation at increased computational cost.
These increases are most important for tori with nearly resonant rotation vectors, high anisotropy, or a wide Fourier spectrum.
For smooth tori with good Diophantine constants, increases of these parameters are not important.
The parameters above give a reasonable balance between speed and accuracy, where the most expensive tori take approximately a 60 seconds to compute.

There are five quantities that we use to analyze the performance of the algorithm
\begin{itemize}
    \item $R_{\mathrm{RRE}}$: The least-squares residual of the Birkhoff RRE procedure \eqref{eq:RRE}. 
    This is an indicator of integrability and a measure of the accuracy of the frequencies $\bm \Omega$.
    \item $L_{\mathrm{MAP}}$: A quantity proportional to the log posterior of the rotation vector MAP problem \eqref{eq:Lmap}.
    This is a measure of confidence in the rotation vector.
    \item $\Rh$: The invariant torus validation error \eqref{eq:validation-error}.
    This is a (rough) measure of the error of the torus parameterization.
    \item $R_{\mathrm{KAM}}$: To validate our results, we use an \textit{a posteriori} residual based on the KAM residual
    \begin{equation}
    \label{eq:kam-residual}
        R_{\mathrm{KAM}} = \frac{\norm{\bm h \circ \sympmap\circ S - \bm h \circ S \circ \tau_{\bm \omega}}_{L^2}}{\norm{\bm h \circ S}_{L^2}}, \qquad \norm{\bm f}_{L^2}^2 = \int_{\Tbb^d}\norm{\bm f(\bm \theta)}^2 \dif \bm \theta,
    \end{equation}
    where the $L^2$ norm is discretized by a $25 \times 25$ trapezoidal quadrature rule.
    Like $\Rh$, this is a measure of the torus parameterization.
    However, unlike $\Rh$, the torus is sampled uniformly rather than along a trajectory.
    \item $M_{\delta}$: The resonance order of the rotation vector $\bm \omega$ \cite{meiss_birkhoff_2021}. 
    For a given tolerance $\delta > 0$, the resonance order of a rotation vector $\bm \omega$ is defined by
    \begin{equation*}
        M_\delta(\bm \omega) = \min_{\bm k \in \Delta(\bm \omega,\delta)} \norm{\bm k}_1, \qquad \Delta(\bm \omega,\delta) = \{\bm k \mid \Tnorm{\bm \omega \cdot \bm k} \leq \delta\}
    \end{equation*}
    That is, it is a measure of the lowest wavenumber such that the rotation vector is nearly resonant.
    So, one might expect a small $M_{\delta}$ would generally correspond to a poorly approximated invariant torus.
    We use the value of $\delta = 10^{-4}$.
\end{itemize}

\paragraph{Choosing the Trajectory Length} 
For the first test, we investigate the required trajectory length for Birkhoff RRE Convergence.
To do this, we scan the filter length $J$ from $4$ to $2500$, choose $T=2J$, and compute $R_{\mathrm{RRE}}$ for each value of $J$ chosen.
The results are shown in Figure \ref{fig:trajectory-classification} (a), where we find that most trajectories have a residual that falls below $10^{-12}$ by $N=8001$ ($J=2000$). 
This cutoff --- plotted as a magenta line --- is the residual by which we determine a trajectory as integrable.
By this test, we classify $851$ of the $1000$ trajectories as integrable with a length $N=8001$ trajectory.

We can compare this to a relative weighted Birkhoff average residual related to the one defined in \cite{sander2020}.
Let $N$ be an even trajectory length. Then, the residual is defined by
\begin{equation*}
    R_{\mathrm{WBA}} = \frac{\norm{\mathcal{WB}_{N/2}[\bm h](F^{N/2}(x)) - \mathcal{WB}_{N/2}[\bm h](x)}}{\norm{\mathcal{WB}_{N/2}[\bm h](x)}}.
\end{equation*}
A plot of this residual as a function of $N$ is shown in Fig.~\ref{fig:trajectory-classification} (b).
The left side of this panel shows the decay of $R_{\mathrm{WBA}}$ over the same trajectory lengths as $R_{\mathrm{RRE}}$.
We see that for the same trajectory length, weighted Birkhoff averaging gives a visually worse classification than Birkhoff RRE.
In the right side of Fig.~\ref{fig:trajectory-classification} (b), we continue the residual $R_{\mathrm{WBA}}$ in blue lines for trajectory lengths up to $N=2 \times 10^5$. 
The weighted Birkhoff classification rate approximately matches the Birkhoff RRE classification rate in the right panel, indicating that weighted Birkhoff averaging requires approximately an order of magnitude longer trajectories for classification.
This observation helps explain why the parameterization coefficients in Sec.~\ref{sec:parameterization} are poorly obtained by weighted Birkhoff averages alone.

In some cases, it may still be preferable to work with weighted Birkhoff averages over Birkhoff RRE for trajectory classification.
Because the cost of a weighted Birkhoff averages is linear in the trajectory length, the residual can be computed much more quickly than Birkhoff RRE. 
This is balanced against the fact that computing the trajectory can be more expensive than performing the RRE least-squares problem (for the $J=2000$ case, the Birkhoff RRE least-squares problem takes approximately $15$ seconds).

\begin{figure}
    \centering
    \includegraphics[width=1.0\linewidth]{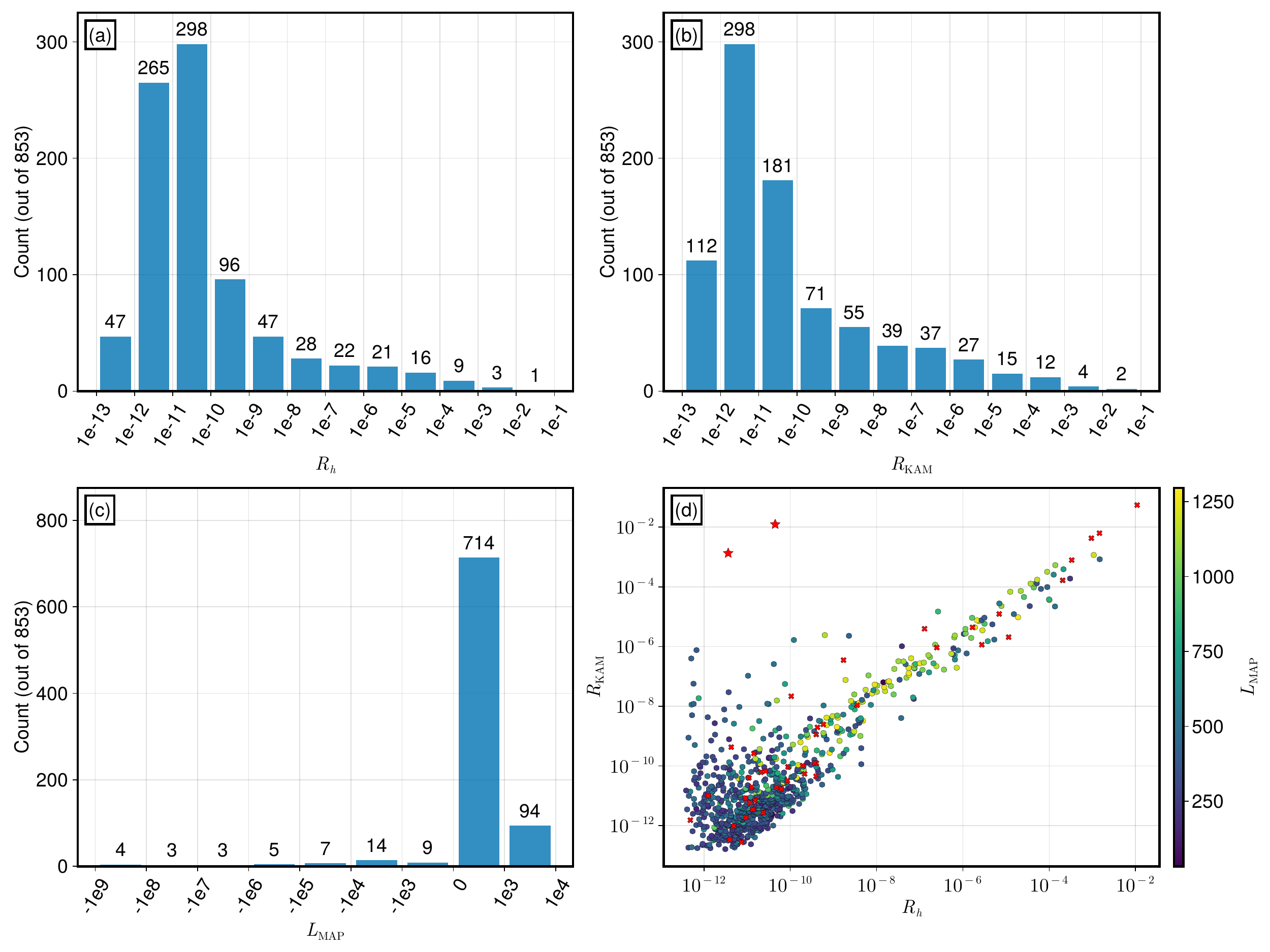}
    \caption{Distribution of the output quantities (a) $\Rh$, (b) $R_{\mathrm{KAM}}$, and (c) $L_{\mathrm{MAP}}$ for the $851$ classified trajectories.
    (d) The validation residual $\Rh$ vs the KAM residual $R_{\mathrm{KAM}}$.
    For $L_{\mathrm{MAP}}>0$, the points are colored according to the value of $L_{\mathrm{MAP}}$.
    For $L_{\mathrm{MAP}}<0$, all but two points are plotted as red `x's.
    The final two are plotted as red stars, indicating the corresponding two the two difficult to approximate trajectories plotted in Fig.~\ref{fig:resonance-lines}. }
    \label{fig:resid-tradeoff}
\end{figure}
\paragraph{Invariant Torus Residuals}
For the $851$ trajectories that are classified as tori by the $J=2000$ case, we apply the methods within Secs.~\ref{subsec:initial-rotation}, \ref{subsec:refine-rotation}, and \ref{subsec:fourier-adaptive} to compute the rotation vector and parameterization.
From this process, we obtain the values of $L_{\mathrm{MAP}}$ and $\Rh$ for each trajectory, and we compute the value of $R_{\mathrm{KAM}}$ as described.

In Fig.~\ref{fig:resid-tradeoff} (a) and (b), we report distributions of $\Rh$ and $R_{\mathrm{KAM}}$ for the trajectories considered.
By both metrics, we find that we have accurately computed parameterizations of most of the trajectories; only six trajectories have $R_{\mathrm{KAM}} \geq 10^{-3}$.
We plot $R_{\mathrm{KAM}}$ against $\Rh$ in Fig.~\ref{fig:resid-tradeoff} (d), and find that there is a strong correlation between the two metrics, indicating that the validation error $\Rh$ can be used as a proxy for the true torus residual in the typical case (two exceptions are plotted as red stars, see the discussion of Fig.~\ref{fig:resonance-lines} below). 

The distribution of $L_{\mathrm{MAP}}$ is given in Fig.~\ref{fig:resid-tradeoff} (c).
We see that most trajectories have $L_{\mathrm{MAP}} > 0$, while there are a few values the negative tail. 
A rough interpretation of this is that, for most trajectories, the assumed Bayesian model well explains the observed data.
We note, however, the connection between this observation and the more fundamental question of whether we determined the correct rotation vector is more complicated.
Typically, the rotation vector is still correct even when $L_{\mathrm{MAP}}$ is negative.
One main reason for this is that islands tend to have negative values of $L_{\mathrm{MAP}}$.
In this case, the posterior could likely be improved by directly working with the map $\sympmap^p$, where $p$ is the period of the island, as this reduces the problem to the invariant torus case.
A second reason is that anisotropic tori can have values of $\bm k$ that are not explained by the value of $P$ chosen in Sec.~\ref{subsec:initial-rotation}.
Of course, this can be improved by choosing $P$ larger.
The final reason is that very smooth tori actually have smaller values of $L_{\mathrm{MAP}}$ in general.

To expand on this final reason, we turn back to Fig.~\ref{fig:resid-tradeoff} (d).
The points where $L_{\mathrm{MAP}} > 0$ are colored according to their values, while values of $L_{\mathrm{MAP}} < 0$ are plotted in red. 
Perhaps against intuition, we find that $L_{\mathrm{MAP}}$ decreases as the torus parameterization becomes more accurate.
The reason is that the more accurately parameterized tori tend to be the smoother ones, and the error in the Birkhoff RRE frequency computation decreases as Fourier modes become more prominent.
For an extreme example of why this might happen, consider a 1D torus of the form $\bm h \circ S(\theta) = (\cos \theta,\sin \theta)$.
Because only the $k = 1$ wavenumber appears in the signal, one can show that a correctly chosen $J=1$ filter gives a residual $R_{\mathrm{RRE}}=0$, regardless of the rotation number. 
When $J$ increases beyond this, nothing remains in the signal that could reveal the rest of the frequencies, so the rest of $\bm \Omega$ would be random.
This means that there is more uncertainty in $\bm \Omega$ for invariant tori that are well described by few Fourier modes, leading to lower values of $L_{\mathrm{MAP}}$. 
Hence, we find that $\Rh$ and $R_{\mathrm{KAM}}$ are much more useful than $L_{\mathrm{MAP}}$ to assess the quality of the final torus.

\begin{figure}
    \centering
    \includegraphics[width=1.0\linewidth]{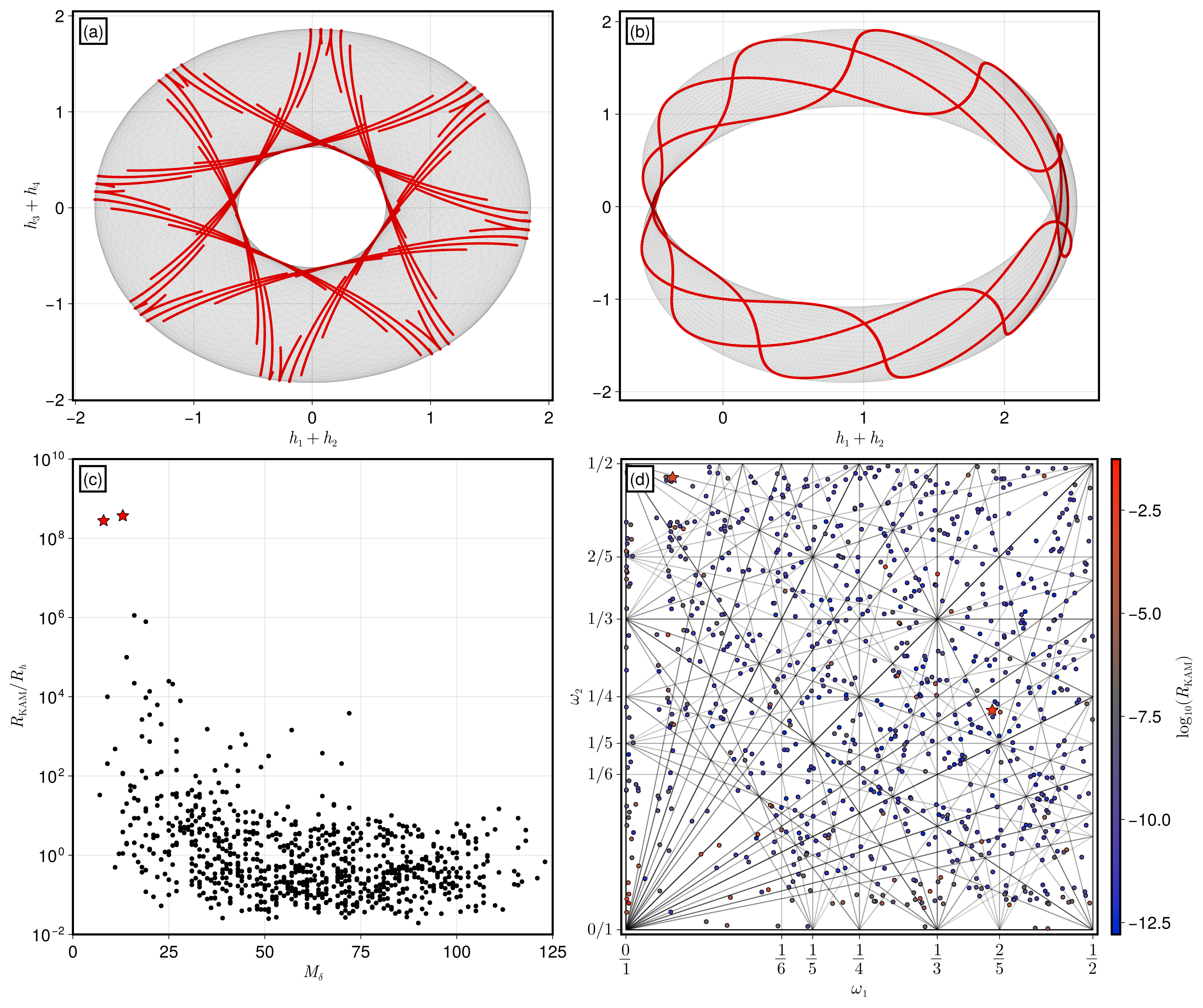}
    \caption{(a-b) Two examples of poorly approximated invariant tori projected onto the $(h_1+h_2, h_3+h_4)$ plane, where the length $N=8001$ trajectories are in red and the computed invariant tori are in gray. 
    (c) A scatter plot of the parameterization error ``surprise'' $R_{\mathrm{KAM}}/\Rh$ against the resonance order $M_\delta$ for each trajectory.
    (d) A scatter plot of the computed rotation vectors $\bm \omega$ over the resonance lines $\bm \theta \cdot \bm k = 0 \mod 1$ for $\norm{\bm k} \leq 6$, colored by the logarithm of the KAM residual.}
    \label{fig:resonance-lines}
\end{figure}

\paragraph{Poorly Approximated Tori and Small Denominators}
We have observed that $\Rh$ tends to be a good proxy for $R_{\mathrm{KAM}}$, but there were two outlier invariant tori (in red stars in Fig.~\ref{fig:resid-tradeoff} (d)).
This is difficult to reconcile, as small $\Rh$ indicates the torus is predictive of the trajectory, but clearly $R_{\mathrm{KAM}}$ is a more honest measure of the torus.
To investigate the mode of failure for the outliers,  in Fig.~\ref{fig:resonance-lines} (a-b) we plot projections of the trajectories in red and the computed invariant tori in black. 
We observe that there are two primary ways in which $R_{\mathrm{KAM}}$ loses accuracy.
For the example Fig.~\ref{fig:resonance-lines} (a), the loss of accuracy is due to a nearly resonant rotation vector of $\bm \omega = [0.39219678, 0.23531362]$, as $\bm \omega \cdot [3,-5] = 2.225 \times 10^{-5} \mod 1$.
The qualitative effect is that the trajectory is slow to explore the surface of the torus, leading to a high condition number for the least-squares parameterization.
For the example in Fig.~\ref{fig:resonance-lines} (b), the loss of accuracy is due to a misidentification of the rotation vector as a consequence of the extreme anisotropy of the torus. 
The computed rotation vector is $\bm \omega = [0.04995223, 0.48501443]$ has a near resonance at $\bm \omega\cdot[3,10]=1.015 \times 10^{-6} \mod 1$, which may indicate that only rotation in the ``long'' direction is being accounted for.
Visually, the trajectory appears as almost periodic.

In both cases, the trajectories could be better approximated by increasing the resolution in the relevant parameters.
For nearly resonant rotation vectors, simply using a longer trajectory for interpolation is sufficient.
For a quicker computation, this trajectory could be subsampled (i.e.~use the map $F^n$ for some integer $n > 1$), as that trajectory would more uniformly explore the torus without increasing the cost of downstream computations.
For highly anisotropic tori, one can increase $J$, $T$, and $\JMAP$.
The use of larger $\JMAP$ gives the algorithm the opportunity to find two rotation vectors that form a valid basis, so that the trajectory can be correctly interpolated.
We note that the higher frequencies are lower accuracy, however, and so there may be a limitation on the current algorithm without also increasing the floating point precision.
As such, the higher anisotropy likely requires further algorithmic improvements for reliable practical computations.

To investigate the issue of near-resonance in a more systematic way, in Fig.~\ref{fig:resonance-lines} (c), we plot the ``surprise'' $R_{\mathrm{KAM}}/\Rh$ against the resonance order $M_\delta$. 
We find that small values of $M_\delta$ are correlated with a mismatch between the validation residual and the KAM residual.
Again, the intuition is likely due to nearly resonant rotation vectors corresponding to worse condition numbers for the parameterization.
In Fig.~\ref{fig:resonance-lines} (d), we plot the rotation vector $\bm \omega$ of each torus, colored by $\log_{10}(R_{\mathrm{KAM}})$.
Behind the rotation vectors, we plot the resonance lines $\bm \theta \cdot \bm k = 0 \mod 1$ for $\abs{\bm k} \leq 6$.
We observe that torus appearing in Fig.~\ref{fig:resonance-lines} (a) lies on the $[5,3]$ resonance line. 
More generally, we find that the KAM residual is typically better within the bulk, while the residual is worse near $\omega_1 = 0$. 
This is likely due to the trend that anisotropic tori and island chains have small $\omega_1$, due to them being near bifurcation points.
Anisotropic tori are hard to approximate because the more complicated Fourier parameterization, while islands chains are complicated because the points per island is divided by the island period.

\begin{figure}
    \centering
    \includegraphics[width=1.0\linewidth]{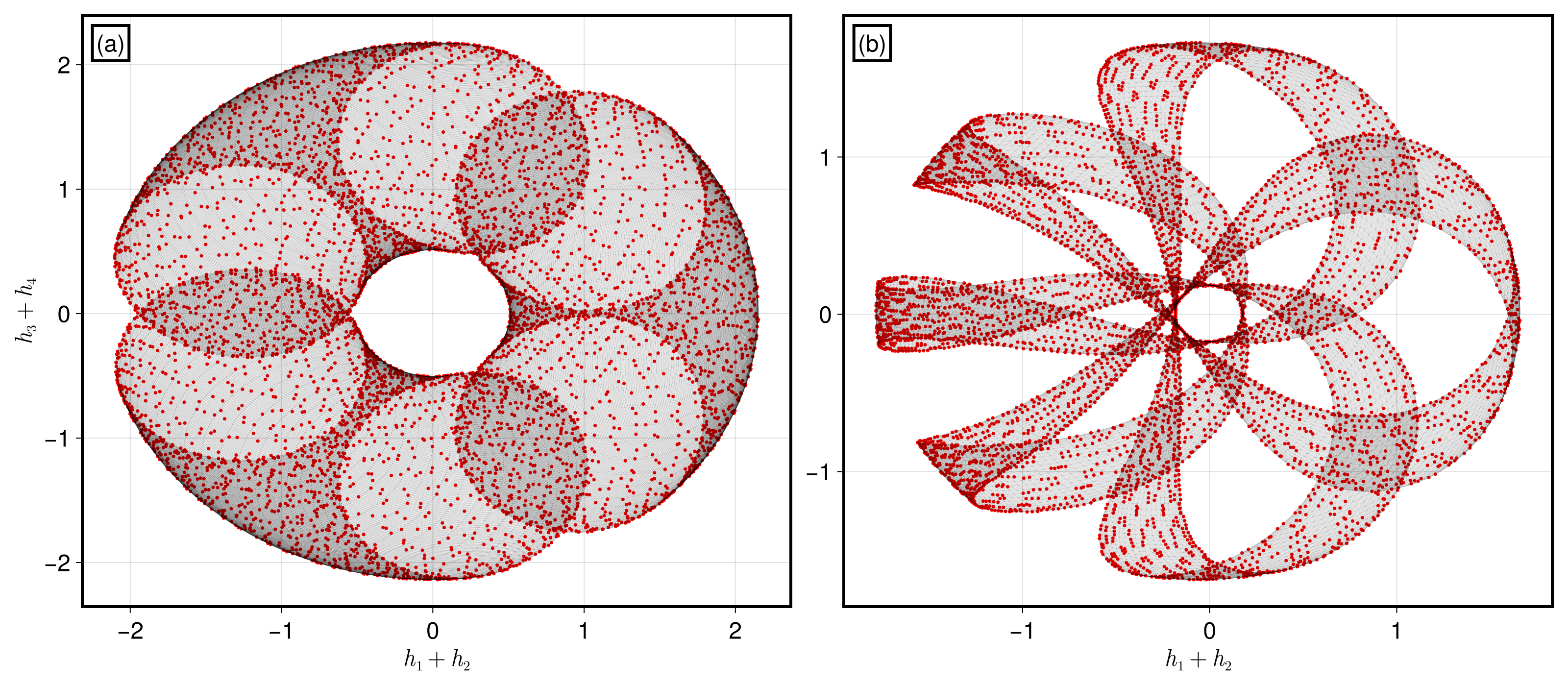}
    \caption{Projections of two successfully computed island trajectories into the $(h_1+h_2, h_3+h_4)$ plane.}
    \label{fig:islands}
\end{figure}

\paragraph{Islands}
We plot two successful parameterizations of island chains in Fig.~\ref{fig:islands} for the coupled standard map.
Fig.~\ref{fig:islands} (a) is a period $p=3$ island and (b) is a period $p=2$ island.
The island structure in (a) resembles a cut donut, while in (b) the islands are significantly more anisotropic.

\subsection{The Earth-Moon restricted three body problem}
\label{subsec:cislunar}
\begin{table}
    \centering
    \begin{tabular}{c|c||c|c|c|c|c}
         Orbit & $N$ & Traj. & RRE LS & RRE e-val & MAP and KZ & Param.  \\ \hline
         $L_4$ Librational & 8335 & 178.7 & 31.8 & 35.3 & 9.5 & 29.0 \\
         Western low prograde & 3335 & 128.9 & 1.4 & 2.4 & 8.0 & 6.8 \\
         Distant retrograde & 3335 & 142.3 & 1.2 & 2.3 & 0.4 & 3.7
    \end{tabular}
    \caption{Timings of each step of computing 3D tori in the earth-moon ER3BP. `Traj.' is the time to compute the trajectory of length $N$, `RRE LS' is the time for the least-squares problem \eqref{eq:RRE}, `RRE e-val' is the time for the eigenvalue problem to obtain frequencies from the RRE filter, `MAP and KZ' is the combined time of the MAP estimation and KZ step (dominated by the MAP step), and `Param.' is the time to compute the adaptive parameterization from Sec.~\ref{subsec:fourier-adaptive}.
    All times are in seconds.}
    \label{tab:3D-timings}
\end{table}

\begin{figure}
    \centering
    \includegraphics[width=1.0\linewidth]{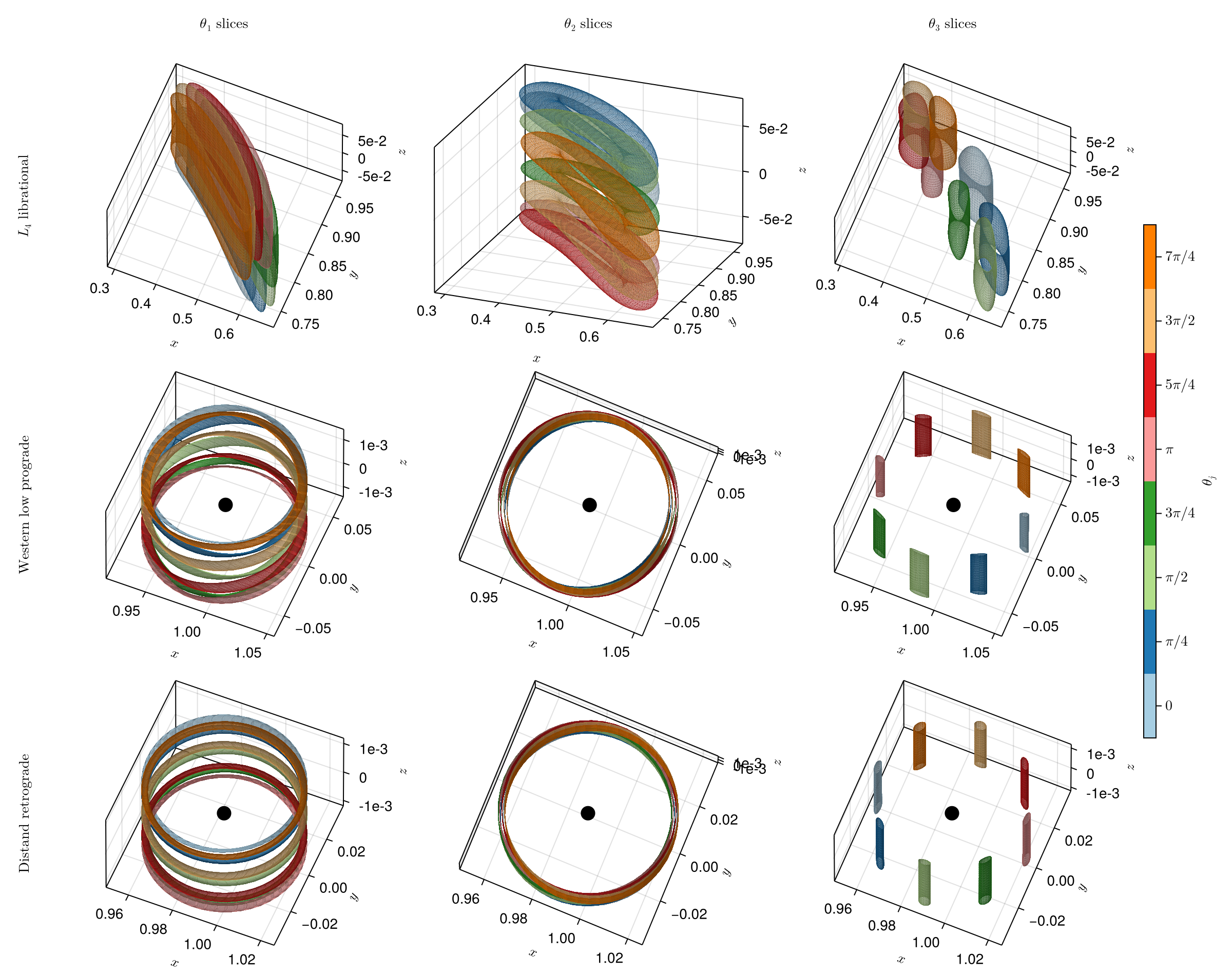}
    \caption{Visualizations of three $d=3$ invariant tori of the earth-moon ER3BP.
    In each column, the three-dimensional torus is ``sliced'' in each of the three angular coordinates. 
    For instance, in the first column, the tori slices $(\bm h \circ S)^{(1,j)}(\theta_2,\theta_3) = (\bm h\circ S)(\pi j / 4, \theta_2, \theta_3)$ for $0\leq j \leq 7$. }
    \label{fig:3D-tori}
\end{figure}

To test our method on 3D invariant tori, we consider the ER3BP with out-of-plane motion.
The differential equations are the same as those presented in Eq. \ref{eq:rer3bp} but with the addition of a third out-of-plane coordinate $\bm \zeta$ so that $\bm\xi=(\xi,\eta,\zeta),\bm p=(p_\xi,p_\eta,p_\zeta)$. 
Accompanying modifications to the Hamiltonian are made accordingly \cite{koon2000dynamical}
\begin{equation*}
    H(\bm \xi, \bm p, f) = \frac{1}{2}\left(\Vert\bm\xi\Vert^2+\Vert\bm p\Vert^2\right) + \eta p_\xi - \xi p_\eta - \phi(\bm \xi,f) \quad \rho_1^2 = (\xi-\mu)^2 + \eta^2+\zeta^2, \quad \rho_2^2 = (\xi - (\mu-1))^2 + \eta^2+\zeta^2.
\end{equation*}
As before, the map $\sympmap$ is obtained by integrating the ODEs with the SRK3 integrator with $1000$ symplectic time steps per lunar orbit.
The observable $\bm h$ is chosen to be the identity. 
We choose the parameters to match the earth-moon orbit with $\mu = 1.2151\times 10^{-2}$ and $\epsilon = 0.0549$.
Because of the increased dimension, the system now supports invariant tori with coupled out-of-plane motion.
We initialize three tori at the locations 
\begin{gather*}
    [\bm \xi_{L_4}, \dot{\bm\xi}_{\mathrm{L_4}}] = \left[\begin{pmatrix} 0.48885 \\ 0.868025 \\ 0.05\end{pmatrix} , \begin{pmatrix}
        -0.002 \\ 0.002 \\0.05
    \end{pmatrix}\right], \quad 
    [\bm \xi_{\mathrm{WLP}}, \bm \dot{\bm\xi}_{\mathrm{WLP}}] = \left[\begin{pmatrix} 1.04254 \\ 0 \\ 0.001 \end{pmatrix} , \begin{pmatrix}
        0 \\ 0.43117 \\ 0.001
    \end{pmatrix}\right], \\
    [\bm \xi_{\mathrm{DR}}, \bm \dot{\bm\xi}_{\mathrm{DR}}] = \left[\begin{pmatrix} 0.95561 \\ 0.001 \\ 0.001 \end{pmatrix} , \begin{pmatrix}
        0.001 \\ 0.64088 \\ 0.001
    \end{pmatrix}\right].
\end{gather*}
where $\dot{\bm\xi}$ is the pulsating rotating frame velocity related to the momenta by
\begin{equation*}
    \bm p =\dot{\bm\xi}+(-\eta,\xi,0)^T.
\end{equation*}
The out-of-plane initial data are perturbations on a $L_4$ librational orbit (similar to  Sec.~\ref{sec:parameterization}), a western low prograde orbit, and a distant retrograde orbit \cite{henon1969numerical, ming2009exploration}.
The initial points were chosen by perturbing orbits catalogued in the NASA JPL Solar System Group Horizons Three-Body Periodic Orbit tool \cite{jpl}.

To parameterize the algorithm, we use (cf.~Sec.~\ref{subsec:sm}) 
\begin{itemize}
    \item$L_4$ \textbf{librational:} $J=2500$, $T=3334$, $N=8335$, $J_0=20$,
    \item\textbf{Western low prograde:} $J=1000$, $T=1334$, $N=3335$, $J_0 = 18$,
    \item\textbf{Distant retrograde:} $J=1000$, $T=1334$, $N=3335$, $J_0 = 8$,
\end{itemize}
and otherwise we keep the same values as the previous section. 
The reported values of $J_0$ for western low prograde and distant retrograde were the minimum values such that $H_{\JMAP} > 10^{-3} H_1$, as described in Sec.~\ref{subsec:initial-rotation}.
The output quantities of interest are:
\begin{itemize}
    \item$L_4$ \textbf{librational:} 
    $\bm \omega = (0.04701, 0.00018, 0.30039)$, $\bm K = (15, 7, 19)$, $R_{\mathrm{RRE}} = 6.28 \times 10^{-11}$, $L_{\mathrm{MAP}} = 938$, $R_h = 1.44 \times 10^{-6}$, $R_{\mathrm{KAM}} = 3.63 \times 10^{-6}$,
    \item\textbf{Western low prograde:} $\bm \omega = (0.03310, 0.47110, 0.11892)$, $\bm K = (7, 11, 25)$, $R_{\mathrm{RRE}} = 8.94 \times 10^{-14}$, $L_{\mathrm{MAP}} = 838$, $R_h = 5.25 \times 10^{-8}$, $R_{\mathrm{KAM}} = 5.70 \times 10^{-6}$,
    \item\textbf{Distant retrograde:} $\bm \omega = (0.36620, 0.16519, 0.40482)$, $\bm K = (7, 13, 21)$, $R_{\mathrm{RRE}} = 7.23 \times 10^{-15}$, $L_{\mathrm{MAP}} = 382$, $R_h = 1.12 \times 10^{-7}$, $R_{\mathrm{KAM}} = 1.63 \times 10^{-5}$.
\end{itemize}
The value of $N$ in the distant retrograde example was chosen to be larger than was necessary to classify the torus, as $\bm \omega$ for appears to be nearly resonant for this trajectory (an issue that we expect becomes worse with dimension).

Each invariant torus is plotted in a row of Fig.~\ref{fig:3D-tori}.
The tori are visualized by plotting ``slices'' of the 3D torus, defined by fixing one of the coordinates $\theta_j$ for $\theta_j =2\pi n/4$ for $0\leq n \leq 7$. 
Each column of Fig.~\ref{fig:3D-tori}, the slicing coordinate is changed.
We see that the third column, slicing in $\theta_3$, aligns with the intuitive notion of the "long way" around the torus, as we might expect from the KZ reduced basis.
However, this is only a rough intuition, as the loops are actually computed in the full six dimensional space.
This is observed through the $\theta_1$ and $\theta_2$ directions, where the shortest way around the torus is more difficult to intuit.

In Table \ref{tab:3D-timings}, we show the wall clock time for each of the steps to compute the tori for each example.
In each case, the time to compute the trajectory was the primary cost.
Second to this, the Birkhoff RRE least-squares problem (computing $\bm c$) and eigenvalue problem (computing $\bm \Omega$ and $\bm H$) scale the worst with the trajectory length, as both use full dense linear algebra with dimensions proportional to the trajectory length.
The parameterization least-squares problem also scales with trajectory length, but the fixed value of $K_{\mathrm{max}}$ improves the performance. 
Finally, because $J_0$ is uncoupled to the trajectory length, the Bayesian MAP problem can dominate the cost at low trajectory length, while being relatively cheap for long trajectories.
In sum, each part of the algorithm can potentially be rate limiting, depending on the context.
\section{Conclusion}
\label{sec:conclusions}
In this paper, we introduced a method that efficiently computes invariant tori from single trajectory of a symplectic map. 
The first step is Birkhoff RRE \eqref{eq:RRE}, which is an extrapolation method that both classifies the torus as integrable or not, and in the case of integrability returns the frequencies and magnitudes of the most prominent Fourier modes.
This frequency content is used by the Bayesian MAP problem \eqref{eq:map-explicit} to compute the valid rotation vector of the invariant torus.
Then, the valid rotation vector is transformed to an optimal one via a KZ basis reduction. 
Finally, using the optimal rotation vector, we perform the optimization \ref{eq:torus-residual} for a parameterization of the invariant torus.

Through the standard map example, we demonstrated that this method can be used to compute large numbers of accurate invariant tori from trajectories of length less than $10000$. 
Most trajectories were computed to a KAM parameterization residual of less than $10^{-4}$. 
This may be enough accuracy for many applications, and if it isn't, the tori could be used to initialize a Newton iteration for a more accurate torus.
We found that the rare tori that fail to be approximated are either highly filamentary or they have nearly resonant rotation vectors.

Looking forward, there are many ways the presented algorithm could be improved. 
From a performance perspective, the four limiting steps are (i) the Birkhoff RRE least-squares problem, (ii) the Birkhoff RRE eigenvalue problem, (iii) the MAP estimation, and (iv) the parameterization optimization. 
Steps (i), (ii), and (iv) could all be improved by applying more efficient numerical linear algebra techniques tailored to the problem structure. 
The matrix for (i) is Hankel, so an iterative least-squares algorithm could be performed with a FFT based algorithm for a worst-case cost of $\Ocal(N^2 \log N)$. 
However, initial tests indicate iterative methods would likely need to be preconditioned to achieve the accuracies necessary for computing the rotation vector.
For (ii), we currently solve the full eigenvalue problem for every root of the filter polynomial, but it is possible that a more efficient algorithm would only compute the few most important frequencies for determining the rotation vector.
For (iv), speedups could potentially be obtained by taking advantage of the discrete Fourier transform-like structure of the linear algebra.
Improvements in step (iii) would need to have a more discrete character, but could potentially be used to find rotation vectors with fewer frequencies of the signal.

From the perspective of theory, a proof of the convergence of the Birkhoff RRE roots to the frequencies $\bm \omega \cdot \bm k$ is still an open question. 
However, the accuracy of the computed parameterizations in this paper supports the conjecture that the roots do converge.
A convergence theorem would also likely require an estimate for the rate of convergence of the roots, which could be incorporated into the MAP problem, and could potentially be pushed all the way through to an error estimate of the final computed torus.

\paragraph{Acknowledgments}
This material is based on work supported by the U.S. Department of Energy, Office of Science, Office of Advanced Scientific Computing Research, as a part of the Mathematical Multifaceted Integrated Capability Centers program, under Award No. DE-SC0023164.

We also thank Amelia Chambliss and Elizabeth Paul.

\appendix
\section{QR Updates for Adaptively Computing Fourier Coefficients}
\label{app:qr-updates}
To improve the efficiency of the adaptive least-squares algorithm in Sec.~\ref{subsec:fourier-adaptive}, we use QR factorization updates \cite{golub_matrix_2013} on the least-squares system \eqref{eq:Fourier-lsqr} at each step.
In particular, in each step of the adaptive algorithm, we need to solve a least-squares problem $A_j(\bm K) X = B$ in each candidate direction $j$, where the matrix $A_j(\bm K)$ consists of columns $\bm \Lambda_{\bm k}$ for $-\bm K - \bm e_j \preceq \bm k \preceq \bm K + \bm e_j$ (see Eq.~\ref{eq:A-Rh}).
As such, we want to avoid recomputing the QR decomposition $A_j(\bm K) = Q_j(\bm K) R_j(\bm K)$ as $\bm K$ increases.

To do so, consider we take a step in the $j'$ direction, having already computed the decompositions $A_j(\bm K) = Q_j(\bm K) R_j(\bm K)$ for each $j$.
For the next step, we need to compute the QR decomposition of $A_j(\bm K + \bm e_{j'}) = Q_j(\bm K + \bm e_{j'})R_j(\bm K+\bm e_{j'})$.
To do this, we notice that the columns of $A_j(\bm K)$ are all contained within $A_j(\bm K + \bm e_{j'})$, so we can write
\begin{align*}
    A_j(\bm K + \bm e_{j'}) &= \begin{pmatrix}
        A_j(\bm K) & A_{j}^{(j')}(\bm K)
    \end{pmatrix} \\
    &= \begin{pmatrix}
        Q_j(\bm K) & Q_{j}^{j'}(\bm K)
    \end{pmatrix}
    \begin{pmatrix}
        R_j(\bm K) & R_j^{(1j')}(\bm K) \\
        0                  & R_j^{(2j')}(\bm K)
    \end{pmatrix},
\end{align*}
where $A_j^{(j')}(\bm K)$ contains the columns in $A_j(\bm K + \bm e_{j'})$ but not in $A_j(\bm K)$, namely $\bm \Lambda_{\bm k}$ for $\bm k$ satisfying $-(\bm K + \bm e_j + \bm e_{j'}) \preceq \bm k \preceq \bm K + \bm e_j +\bm e_{j'}$ and $\abs{k_{j'}} = (\bm K + \bm e_j)_{j'} + 1$.
Then, the unknown QR factors can be computed by the update formulas
\begin{align*}
    R_j^{(1j')} &= Q_j^* A_j^{(j')}, \\
    Q_j^{(j')} R_j^{(2j')} &=A_j^{(j')} - Q_j R_j^{(1j')},
\end{align*}
where we suppressed the argument of $\bm K$ for all matrices and the second line is computed by a standard QR factorization.
We note that the asymptotic cost to compute a QR factorization by updates is equal to the cost of the QR factorization, leading to the complexity estimates in Sec.~\ref{subsec:fourier-adaptive}.

\section{Solving the RRE Least-Squares System}
\label{app:solving-RRE}
In order to solve the constrained least-square system \eqref{eq:RRE}, we perform a minor update on the original method \cite{ruth_finding_2024}. 
The main idea is for the palindromic and sum-to-one constraints to be satisfied to near machine precision, one must perform an orthogonal transformation to obtain a solution.
Define the discrete cosine transform matrix (ignoring the constant mode)
\begin{equation*}
    Q = \begin{pmatrix}
        | & & |\\
        \bm q_1 & \dots & \bm q_{J} \\
        | & & |
    \end{pmatrix}, \qquad (\bm q_k)_j = \sqrt{\frac{2}{2J+1}} \cos\left(2\pi k \frac{j}{2J+1} \right), \qquad -J\leq j\leq J.
\end{equation*}
If we let $\bm c = (2J+1)^{-1}\bm 1 + Q \bm c'$ where $\bm 1$ is the vector of all ones, the fact that $(\bm q_k)_{-j} = (\bm q_k)_j$ clearly implies that $\bm c$ is palindromic.
Moreover, the orthogonality of the discrete cosine transform (i.e.~$\bm 1 \cdot \bm q_k = 0$) implies the sum-to-one constraint is satisfied.
So, we can substitute this representation into the least-squares problem to find
\begin{equation*}
    R_{\mathrm{RRE}}^2 = \min_{\bm c' \in \Rbb^{J}} \frac{1}{T} \norm{U Q\bm c' + \frac{1}{2J+1} U \bm 1}^2.
\end{equation*}
We solve this by a direct QR factorization on $U Q$.

\bibliographystyle{unsrt}
\bibliography{bibliography}

\end{document}